\documentclass[a4paper]{article}
\usepackage{amsmath}
    \usepackage{amsxtra}
    \usepackage{amstext}
    \usepackage{amssymb}
    \usepackage{latexsym}

\input xy

\xyoption{all}

\xyoption{poly}

\newcommand{\nin}{\mbox{$\in \!\!\!\!\!/$}}

\newtheorem{prop}{Proposition}[section]
\newtheorem{thm}{Theorem}[section]
\newtheorem{cor}{Corollary}[section]
\newtheorem{lem}{Lemma}[section]

\newtheorem{cl}{Claim}[section]

\newtheorem{deff}{Definition}[section]

\begin{document}

{\begin{center}\large \bf NORMAL PROJECTIVITY OF COMPLETE\\
SYMMETRIC VARIETIES
\end{center}}

{\begin{center} \sc ALESSANDRO RUZZI
\end{center}}

\begin{abstract}
Chiriv\`{\i} and Maffei \cite{CM II} have proved that the
multiplication of sections of any two ample spherical line bundles
on the wonderful symmetric variety $X=\overline{G/H}$ is surjective.
We have proved two criterions that allows ourselves to reduce the
same problem on a (smooth) complete symmetric variety to the
corresponding problem on the complete toric variety (respectively to
the open toric variety). We have also studied in details some family
of complete symmetric varieties, in particular those of rank 2.
\end{abstract}


In this work we will study the projective normality of  complete
symmetric varieties. Let $\overline{G}$ be an adjoint semisimple
group over $\textbf{C}$ and let $\theta$ be an involution of
$\overline{G}$. We define $\overline{H}$ as the subgroup of the
elements fixed by $\theta$ and we will say that
$\overline{G}/\overline{H}$ is a homogeneous symmetric variety. De
Concini and Procesi \cite{CSV I}  have defined a wonderful
completion of $\overline{G}/\overline{H}$ and this is the unique
wonderful completion of $\overline{G}/\overline{H}$. In this work we
will define a complete symmetric variety as an
$\overline{G}$-variety with a dense open orbit isomorphic to
$\overline{G}/\overline{H}$ and a $\overline{G}$-equivariant map
$Y\rightarrow X$ extending the identity of
$\overline{G}/\overline{H}$. They are be classified by De Concini
and Procesi \cite{CSV II}. Indeed they showed that there is an
equivalence of categories between the category of complete symmetric
varieties  and the category of  toric varieties over an affine space
$\textbf{A}^{l}$ considered as a $(\textbf{C}^{*})^{l}$ variety in
the obvious way, where $l$ is the rank of $G/H$. Moreover there is a
one-to-one correspondence between the completions $Y$ of
$\overline{G}/\overline{H}$ which lie over $X$ and the elements of a
special class of  complete toric varieties. One can show that the
complete toric variety $Z^{c}$ corresponding to a complete symmetric
variety $Y$ is a subvariety of $Y$ and the open toric variety $Z$
corresponding to $Y$ is an open subvariety of $Z^{c}$.

In this work, unless explicitly stated, we shall always assume that the complete symmetric variety $Y$ is
smooth. Recall that by $\cite{CSV II}$ it then follows that: 1) any orbit closure in $Y$ is also smooth; 2) the
associated toric varieties $Z$ and $Z^{c}$ are both smooth.

Our work consist of two parts. In the first one we will reduce the
study of the projective normality of the complete symmetric
varieties to the study of projective normality of the corresponding
complete toric varieties. First of all we will prove that a complete
symmetric variety is projective if and only if the corresponding
complete toric variety is projective. Moreover they are projective
if and only if the associated open toric variety is
quasi-projective.

Chiriv\`{\i} and Maffei \cite{CM II} have proved the following
result which easily implies the projective normality of $X$ with
respect to any projective embedding by a complete linear system.

\begin{thm}\label{CM}
Let $L_{1}$ and $L_{2}$ be any two line bundles generated by global sections on the wonderful complete symmetric
variety $X$. Then the product of sections
\[H^{0}(X,L_{1})\otimes H^{0}(X,L_{1})\longrightarrow H^{0}(X,L_{1}\otimes L_{2})\]
is surjective.
\end{thm}

We will try to generalize this results to any complete symmetric
variety. First we will prove that the surjectivity of the product of
sections of two ample line bundles on a complete symmetric variety
is equivalent to the surjectivity of the product of  sections of the
restrictions of the line bundles to the corresponding complete toric
variety. Thus we will have reduced the problem to a problem on toric
varieties. Unfortunately, it is very difficult to verify the
surjectivity of the product of sections of any two ample line
bundles on a generic complete toric variety. However, we can
simplify the problem for the special class of complete toric
varieties which we are considering. Indeed we will prove that the
surjectivity of the product of sections of two ample line bundle on
$Z^{c}$, say $L_{1}$ and $L_{2}$, is equivalent to the surjectivity
of the product of sections of the restrictions of the line bundles
to $Z$. This problem is much simpler, because $H^{0}(Z,L_{1}|Z)$ and
$H^{0}(Z,L_{2}|Z)$ are infinite dimensional vector spaces and it is
sufficient to prove that  a suitable finite dimensional subspace of
$H^{0}(Z,(L_{1}\otimes L_{2})|Z)$ is contained in the image of the
product of sections. Indeed we will prove that, given any ample line
bundle $L$ on $Z^{c}$, $H^{0}(Z,L |Z)$ is generated by
$H^{0}(Z^{c},L |Z^{c})$ as an $\mathcal{O}_{Z}(Z)$-module.

In the second part of this work we will study the suriectivity of
the product of sections of an  ample line bundle on a  toric variety
proper over $\textbf{A}^{l}$. If the variety has dimension 2, we
will prove  the product of sections of any two ample line bundles is
surjective. If the dimension of the variety is larger than 2, we
will find a infinite number of varieties proper over
$\textbf{A}^{l}$ such that, for any ample line bundle $L$ on a such
variety $Z$, the product of sections of $L$
\[H^{0}(Z,L)\otimes H^{0}(Z,L)\longrightarrow H^{0}(Z,L\otimes L)\]
is surjective.

We would like to thank C. De Concini for the continuous help.
Moreover we would like to thank W. Fulton for some useful
information.

\section{Introduction}

Let $G$ be a semisimple simply connected algebraic group over $\textbf{C}$ and let $\theta$ be an involution of
$G$. We define $H$ as the normalizer $N_{G}(G^{\theta})$ of the subgroup of  $\theta$-fixpoints. Let
$\overline{G}$ be the adjoint semisimple group associated to $G$ and let $\overline{H}$ be the subgroup of the
elements fixed by the involution induced by $\theta$, then $\overline{G}/\overline{H}$ is isomorphic to $G/H$
through the map induced by the quotient map $G\rightarrow\overline{G}$.

\begin{deff}
We will say that $G/H$ is a homogeneous symmetric variety (of adjoint type).
\end{deff}

We can associate a possibly not reduced root system to $G/H$ (see
\cite{He}). Let $T^{1}$ be a torus of $G$ such that: 1)
$\theta(t)=t^{-1}$ for each $t\in T^{1}$; 2) the dimension $l$ of
$T^{1}$ is maximal. Let $T$ be a maximal torus which contains
$T^{1}$. One can show that $T$ is stabilized by $\theta$, so
$\theta$ induces an involution on $\chi^{*}(T)_{\textbf{R}}$, which
we call again $\theta$. This involution stabilizes the root system
$\phi$ of $G$ and it is orthogonal with respect to the Killing form.
We can choose a Borel subgroup $B$ of $G$ such that the associated
set of positive roots $\phi^{+}$ has the following property: for
each $\alpha\in\phi^{+}$ either $\theta(\alpha)$ is equal to
$\alpha$ or $\theta(\alpha)$ is a negative root. For each root
$\alpha$ we define $\alpha^{s}=\alpha-\theta(\alpha)$. The set
$\{\alpha^{s}\neq0: \alpha\in\phi\}$ is a possibly not reduced root
system of rank $l$ called the restricted root system. We will say
that the not-zero $\alpha^{s}$ are restricted roots and that $l$ is
the rank of $G/H$. The restricted roots generate the
$(-1)$-eigenspace of $\chi^{*}(T)_{\textbf{R}}$.

We now want to describe the lattice of integral weights of $\phi^{s}$. Let $\Lambda$ be the lattice of integral
weights of $\phi$ and let $\Lambda^{+}$ be the set of dominant weights. For each $\alpha\in\Gamma$, let
$\omega_{\alpha}$ be the fundamental weight associated to $\alpha$. Moreover, for each dominant weight
$\lambda$, let $V_{\lambda}$ be the irreducible representation of $G$ of highest weight $\lambda$. We will say
that a dominant weight $\lambda$ is spherical if $V_{\lambda}$ contains a not-zero vector fixed by the Lie
algebra $\mathfrak{h}$ of $H$. Moreover we will say that a weight $\mu$ is special if $\theta(\mu)=-\mu$. Let
$\Omega$ be the lattice generated by the spherical weights and let $\Lambda^{1}$ be the lattice of the special
weights. One can easily show that $2\Lambda^{1}\subset\Omega\subset\Lambda^{1}$. We can describe $\Omega$ more
explicitly. Let $\phi_{0}$ be the set of root fixed by $\theta$ and let $\phi_{1}=\phi-\phi_{0}$. We set
$\Gamma_{0}=\Gamma\cap\phi_{0}$ and $\Gamma_{1}=\Gamma\cap\phi_{1}$ where $\Gamma$ is the basis of $\phi$
associated to $\phi^{+}$. The set $\Gamma^{s}=\{\alpha^{s}: \alpha\in\Gamma_{1}\}$ is a base of $\phi^{s}$.
$\theta$ induces an involution $\overline{\theta}$ of $\Gamma_{1}$ such that, for each $\alpha\in\Gamma_{1}$,
$\theta(\alpha)=-\overline{\theta}(\alpha)-\beta_{\alpha}$ where $\beta_{\alpha}$ is a linear combination with
positive integral coefficients of simple roots fixed by $\theta$. We can order the simple roots
$\alpha_{1},...,\alpha_{l}, \alpha_{l+1},...,\alpha_{l+s},\alpha_{l+s+1},...,\alpha_{m}$ so that $\alpha_{i}$ is
fixed by $\theta$ if and only if $i>l+s$. Moreover we can suppose that $\alpha^{s}_{1},...,\alpha_{l}^{s}$ are
distinct. For each $i=1,...l$ we define $\omega_{i}$ as follows: if $\overline{\theta}(\alpha)=\alpha$ then
$\omega_{i}=\omega_{\alpha_{i}}$, otherwise $\omega_{i}=\omega_{\alpha_{i}}-\theta(\omega_{\alpha_{i}})=
\omega_{\alpha_{i}}-\omega_{\overline{\theta}(\alpha_{i})}$. A weight is special if and only if it is a linear
combination of the weight $\omega_{1},...,\omega_{l}$, so $\{\omega_{1},...,\omega_{l}\}$ is a basis of
$\Lambda^{1}$. Moreover we can use these weights to give the  following explicit description of $\Omega$.

\begin{prop}[Lemma 2.1 and Theorem 2.3 in \cite{CM I}]\label{CM Ithm}
Let \, $\Omega^{+}$ be the set of spherical weights, then
$\Omega\cap\Lambda^{+}=\Omega^{+}$. Moreover
$\Omega=\bigoplus_{i=1}^{l}\textbf{Z}a_{i}\omega_{i}$ where
$a_{i}\in \{1,2\}$ for each $i$. $a_{i}$ is equal to 2 if
$\theta(\alpha_{i})= -\alpha_{i}$, while it is equal to 1  if
$\theta(\alpha_{i})\neq -\alpha_{i}$. In particular $a_{i}=1$ if
$\overline{\theta}(\alpha_{i})\neq \alpha_{i}$. For each $i$ and $j$
we have $<a_{i}\omega_{i},
(\alpha^{s}_{j})^{\vee}>=b_{i}\delta_{i,j}$ where
$(\alpha^{s}_{j})^{\vee}$ is the coroot associated to
$\alpha^{s}_{j}$ and $b_{i}\in\{1,2\}$. $b_{j}=2$ if and only if
$2\alpha^{s}_{j}\in\widetilde{\phi}$. In particular, if
$\widetilde{\phi}$ is reduced then
$a_{1}\omega_{1},...,a_{l}\omega_{l}$ are the fundamental weights
dual to $(\alpha^{s}_{1})^{\vee},..., (\alpha^{s}_{l})^{\vee}$.
\end{prop}

Notice that the proposition implies that the fundamental Weyl
chamber $C^{+}$ of $\widetilde{\phi}$ is the intersection of
$M_{\textbf{R}}$ with the fundamental Weyl chamber $\Lambda^{+}$ of
$\phi$. We will say that a special weight $\sum n_{i}\omega_{i}$ is
regular if $n_{i}>0$ for each $i$. Thus a spherical weight is
regular if and only if it is a strongly dominant weight of the
restricted root system (with respect to the basis $\Gamma^{s}$). The
Weyl group $W^{1}$ of the restricted root system   is called the
restricted Weyl group and it has the following description.

\begin{prop}[See 1.8 in \cite{dCS}]\label{hsv3}
One can identify the restricted Weyl group   $W^{1}$ with the group $\{w\in W:w\cdot
\mathfrak{t}_{1}\subset\mathfrak{t}_{1}\}/W_{0}$, where $W$ is the Weyl group of $\phi$ and $W_{0}$ is the Weyl
group of the root system $\phi_{0}$ (in $\chi^{*}(T^{0})_{\textbf{R}}$).
\end{prop}

De Concini and Procesi have defined a wonderful compactification of $G/H$. Let $\lambda$ be a regular spherical
weight and let $k_{\lambda}$ be a not-zero vector of $V_{\lambda}$ fixed by $\mathfrak{h}$. One can show that
$k_{\lambda}$ is unique up to a not-zero scalar. Let $x_{0}$ be the class of $k_{\lambda}$ in
$\textbf{P}(V_{\lambda})$. De Concini and Procesi have defined the wonderful compactification $X$ of $G/H$ as
the closure of $Gx_{0}$ in $\textbf{P}(V_{\lambda})$.

We now want to give a local description of $X$. Let $T^{0}$ be the
connected component of the subgroup of the invariant of $T$ and let
$S$ be the quotient  of $T^{1}$ by $T^{1}\cap T^{0}$. Choose a basis
of weight vectors and consider the affine open set  $\widetilde{A}$
of $\textbf{P}(V)$ where the coordinate corresponding to the highest
weight $v_{\lambda}$ is not zero. Observe that $\widetilde{A}\cap X$
is $U^{-}$ stable, where $U^{-}$ is the unipotent group associated
to $-\phi_{1}^{+}$, namely
$U^{-}=\prod_{-\alpha\in\phi_{1}^{+}}U_{\alpha}$ as a variety. One
can show that the closure of $T[k]$ in $\widetilde{A}$ is  an affine
space $\textbf{A}^{l}$ with coordinates
$-\alpha^{s}_{1},...,-\alpha^{s}_{l}$. Moreover the map
$\varphi:U^{-}\times \textbf{A}^{l}\rightarrow \widetilde{A}\cap X$
given by $\varphi(g,v)=g\cdot v$ is an isomorphism. For each $i$,
let $X_{-\alpha^{s}_{i}}$ be the divisor of $X$ whose intersection
with $U^{-}\times \textbf{A}^{l}$ is the locus of zeroes of
$-\alpha^{s}_{i}$. De Concini and Procesi have proved that $X$ is
wonderful according to the definition of Luna:

\begin{thm}[Theorem 3.1 in \cite{CSV I}]
\label{CSV I 3.1} Let  $X$ be as before, then:

\begin{enumerate}
\item the stabilizer of $x_{0}$ is $H$;

\item $X$ is a smooth  projective   $G$-variety;

\item $X\backslash (G\cdot [k])$ is a divisor with normal  crossings. It has irreducible components
$X_{-\alpha^{s}_{1}},...,X_{-\alpha^{s}_{l}}$ and they are smooth subvarieties of $X$;

\item the $G$-orbits of $X$ correspond to the subsets of
$\{1,2,...,l\}$,  so that the orbit closures are the intersections $X_{-\alpha^{s}_{i_{1}}}\cap ...\cap
X_{-\alpha^{s}_{i_{k}}}$;

\item there is an unique closed orbit $\bigcap_{i=1}^{l}X_{-\alpha^{s}_{i}}$ and it is
isomorphic to  $G/P$, where $P$ is the parabolic subgroup of $G$ associated to $\Gamma_{0}$, i.e. its Lie
algebra is $\mathfrak{t}\oplus\bigoplus_{\alpha\in\phi_{0}\cup\phi^{+}}\mathfrak{g}_{\alpha}$;

\item $X$ does not depend by $\lambda$.
\end{enumerate}

\end{thm}

De Concini and Procesi have also classified the embeddings of $G/H$ over $Y$ i.e. the $G$ varieties $Y$ with a
open orbit isomorphic to $G/H$ and an equivariant map $\pi:Y\rightarrow X$ extending the identity map of $G/H$.
We will say that a such variety is a complete symmetric variety if it is complete.

\begin{thm}
There is an equivalence between the category of complete symmetric varieties   and the category of  toric
varieties proper over $\textbf{A}^{l}$. Given a complete symmetric variety $\pi:Y\rightarrow X$ the
corresponding open toric variety $Z$ is the inverse image of $\textbf{A}^{l}$ in $Y$. Moreover $Y$ is smooth if
and only if  $Z$ is smooth.
\end{thm}

One can show that the inverse image of the open set $U^{-}\times\
\textbf{A}^{l}$ is isomorphic to
$U^{-}\times\pi^{-1}(\textbf{A}^{l})$ in a $U^{-}\cdot T$
equivariant way. Moreover the $G$ orbits of $Y$ are in one-to-one
correspondence to the $S$ orbits of $Z$. The linear map
$\xymatrix{\chi^{*}(S)_{\textbf{R}} \ar@{->}[r]
&\chi^{*}(T^{1})_{\textbf{R}}}$ induced in an obvious way by the
quotient map $\xymatrix{ T^{1}\ar@{->>}[r]&S  }$ is an isomorphism,
while the linear map $\xymatrix{
\chi^{*}(T)_{\textbf{R}}\ar@{->}[r]&\chi^{*}(S)_{\textbf{R}}}$
induced by the canonical injection $\xymatrix{T^{1} \ar@{^{(}->}[r]
&T}$ is surjective. Thus we obtain a surjective map $\xymatrix{
\chi^{*}(T)_{\textbf{R}}\ar@{->>}[r]&\chi^{*}(S)_{\textbf{R}}  }$,
whose restriction to the \mbox{(-1)-eigenspace} is an isomorphism.
This map allow ourselves to identify  $\chi^{*}(S)$ with the lattice
$M$ generated by the restricted roots. We call $N$ the dual group of
$M$, i.e. $N=Hom(M,\textbf{Z})$.

De Concini and Procesi have also proved that there is an one-to-one correspondence between the complete
symmetric varieties and a class of complete toric varieties. The closure $Z_{0}^{c}$ of $\textbf{A}^{l}$ in $X$
is the complete toric variety whose fan is composed by the Weyl chambers and their faces.
\begin{thm}
There is an one to one correspondence between complete symmetric varieties and complete toric $S$-varieties over
$Z_{0}^{c}$ whose fan is $W^{1}$ invariant.
\end{thm}
Given a complete symmetric variety $\pi:Y\rightarrow X$ the corresponding complete  toric variety $Z^{c}$ is the
inverse image of $Z_{0}^{c}$ in $Y$. One can show that $Z^{c}$ is the closure of $Z$ in $Y$. Moreover $Y$ is
smooth if and only if $Z^{c}$ is smooth. {\em In this work, unless explicitly stated, we shall always assume
that the complete symmetric variety $Y$ is smooth.} In this case any orbit closure in $Y$ is also smooth.

We now want to study the line bundles on a complete symmetric variety. We begin with the wonderful case. First
of all, we can identify the Picard group of the wonderful symmetric variety $X$ with a subgroup $\Lambda_{X}$ of
the lattice $\Lambda$ of integral weights of $G$. Indeed De Concini and Procesi have proved that:

\begin{prop}[Proposition 8.1 in \cite{CSV I}]\label{CSV I 8.1}
The map $Pic(X)\rightarrow Pic(G/P)$ induced by the canonical inclusion is injective.
\end{prop}
Remember that we can identify $Pic(G/P)$ with a sublattice of the lattice of weights. $Pic(G/P)\equiv
Pic_{G}(G/P)$ because $G/P$ is a spherical variety. Thus to any linearized line bundle $L\in Pic_{G}(G/P)$ we
associate the opposite $\lambda$ of the character $-\lambda$ with which $T$ acts on the fibre over $P\in G/P$.
Because of the previous proposition, we denote a line bundle  on $X$ with $L_{\lambda}$ if its image is the
weight $\lambda$. Let $\lambda$ be a dominant  weight such that $\textbf{P}(V_{\lambda})$ contains a line $r$
fixed by $H$, for example $\lambda\in\Omega^{+}$. One can show that the map $G/H\ni gH\rightarrow gr$ can be
extended to a morphism $\psi_{\lambda}:X\rightarrow \textbf{P}(V_{\lambda})$. The line bundle
$\psi_{\lambda}^{*}O(1)$ is $L_{\lambda}$. Thus, given a line bundle $L_{\mu}$ on $X$ such that $\mu$ is
dominant, there is a sub-representation of $H^{0}(X,L_{\mu})$ isomorphic to $V_{\mu}^{*}$ and  obtained by
pullback of $H^{0}(\textbf{P}(V_{\mu}),\mathcal{O}(1))$ to $X$. $H^{0}(X,L_{\mu})$ is multiplicity free because
$X$ is a spherical variety, so the previous representation is unique and  we can call it $V_{\mu}^{*}$. Before
to give an explicit description of $Pic(X)$ we need a definition.

\begin{deff} We will say that a root $\alpha\in\Gamma_{1}$ is an exceptional root if
$\overline{\theta}(\alpha)\neq\alpha$ and $<\alpha,\theta(\alpha)>\neq0$.  Moreover we will say that $G/H$ is
exceptional if  there is an exceptional root. We will say that a compactification of $G/H$ is exceptional if
$G/H$ is exceptional.
\end{deff}

\begin{thm}[Theorem 4.8 in \cite{dCS}]
\label{dCS 4.8} $Pic(X)$ is generated by the spherical weights and by the fundamental weights corresponding to
the exceptional roots.
\end{thm}

We now consider the general case following \cite{B}. First of all we
introduce some notations that we will often use.  Let $X$ be the
wonderful complete variety and let $Y$ be the complete symmetric
variety over $X$ associated to a toric variety $Z$ proper over
$\textbf{A}^{l}$. We will call $\Delta$  the fan of $Z$ and
$\Delta^{c}$  the fan of $Z^{c}$. We shall denote  the fan of
$Z_{0}:=\textbf{A}^{l}$ by $\Delta_{0}$ and  the fan of $Z^{c}_{0}$
by $\Delta^{c}_{0}$. Let  $o_{\gamma}$ be the $S$-orbit of $Z$
associated to $\gamma\in\Delta$. We will call $O_{\tau}$  the
$G$-orbit of $Y$ corresponding to $o_{\tau}$. We shall denote by
$Z_{\gamma}$  the stable subvariety of $Z$ associated to
$\gamma\in\Delta$, by $Z^{c}_{\gamma}$  the stable subvariety of
$Z^{c}$ associated to $\gamma\in\Delta^{c}$ and by $Y_{\gamma}$  the
stable subvariety of $Y$ associated to $\gamma\in\Delta$. We set
$\Delta(i)=\{\gamma\in\Delta: \textrm{dim}\,\gamma=i\}$ and
$\Delta^{c}(i)=\{\gamma\in\Delta^{c}: \textrm{dim}\,\gamma=i\}$.

Remember that the closed orbits $O_{\sigma}$ of $Y$ are in one-to-one correspondence with the maximal cones of
the fan $\Delta$ associated to $Z$. Moreover they are all isomorphic to the unique closed orbit of $X$ through
the restriction of the projection, so we can identify $Pic(O_{\sigma})$ with $\Lambda_{X}$ for each
$\sigma\in\Delta(l)$. One can easily show that $Pic(Z)$ is freely generated by the line bundles $O(Z_{\tau})$
where $\tau$ varies in the set  $ \Delta(1)\backslash \Delta_{0}(1)$.  The following theorem gives a complete
description of $Pic(Y)$. Given a cone $\gamma\subset N_{\textbf{R}}$ we will call $\gamma^{\bot}$ the subspace
of $M_{\textbf{R}}$ of the vectors which vanishes on $\gamma$.

\begin{thm}[Theorem 2.4 in \cite{B}] \label{B 2.4}
Let $Y$ be the complete symmetric variety associated to $Z$. Then
\begin{enumerate}
\item The maps $\xymatrix{Z  \ar@{^{(}->}[r]^{i}& Y\ar[r]^{\pi} &X}$ induce the split exact sequence

\[\xymatrix{\ \  0\ar[r] & Pic(X)\ar[r]^{\pi^{*}}& Pic(Y) \ar[r]^{i^{*}} &Pic(Z)\ar[r]&0,}\]

so $Pic(Y)$ is (not canonically) isomorphic to $Pic(X)\oplus Pic(Z)$.

\item A section  of the previous split short exact sequence
is given by sending the free generators $\mathcal{O}(Z_{\tau})$, with $\tau\in \Delta(1)\backslash
\Delta_{0}(1)$, to $\mathcal{O}(Y_{\tau})$. Thus

\[ Pic(Y)= \pi^{*}Pic(X)\oplus \bigoplus_{\tau\in \Delta(1)\backslash \Delta_{0}(1)} \textbf{Z}\,
\mathcal{O}(Y_{\tau}).\]

\item The morphism given by the restriction to the closed orbits

\[c_{1}^{G}:Pic(Y)\rightarrow \prod_{\sigma\in\Delta(l)} Pic(O_{\sigma})\]
is injective and its image can be identified with the lattice

\[\Lambda_{Y}= \{h=(h|\sigma)\in \prod_{\sigma\in\Delta(l)} \Lambda_{X}\subset \prod_{\sigma\in\Delta(l)}
\Lambda:\ h|\sigma-h|\sigma' \in M\cap (\sigma\cap\sigma')^{\bot}\]\[ \forall\ \sigma,\ \sigma'\, \in\Delta(l) .
\}\]
\end{enumerate}
\end{thm}

We will indicate with $L_{h}$ the line bundle whose image is $h$. $Pic(Y)$ is isomorphic to the group of
equivariant line bundles $Pic_{G}(Y)$ because  $G$ is simply connected and $Y$ is complete. Moreover, given a
line bundle $L_{h}$, $-h_{\sigma}$ is the character of the action of $T$ on the fibre over the $T$-stable point
$O_{\sigma}\cap Z$. In a similar way, we define $h^{c}$ as the set $(h|\sigma)$ where $\sigma$ varies in
$\Delta^{c}(l)$ and $-h_{\sigma}$ is the character of the action of $T$ on the fibre over the $T$-stable point
of $Z^{c}$ corresponding to $\sigma$. In some case we can give an useful interpretation of $h$ and $h^{c}$.
Before  doing it, we need to give some definitions. Let $M'$ be a lattice in $M_{\textbf{R}}$ which contains
$M$.
\begin{deff} A real valued function
$h:|\Delta|\rightarrow\textbf{R}$ on the support of $\Delta$ is called a $(\Delta,M')$-linear function if $h$ is
linear on each $\sigma\in \Delta$. Let $h|\sigma$ be the unique linear function which coincide with $h$ on
$\sigma$. We request moreover that $h|\sigma$  belongs to $M'$ and that $h|\sigma_{1}-h|\sigma_{2}$ belongs to
$M$ for each $\sigma$, $\sigma_{1}$ and $\sigma_{2}$ in $\Delta(l)$. Let $SF(\Delta,M')$ be the additive group
of the $(M',\Delta)$-linear functions.
\end{deff}
Sometimes we say that an element  $h$ of $ SF(\Delta,M')$ is a $\Delta$ linear function. We can  think a
$\Delta$-linear function as a function $h:N_{\textbf{R}}\rightarrow \textbf{R}\cup\{-\infty\}$. We say that a
$\Delta$ linear function $h$ is convex if $h(v+v')\geq h(v)+h(v')$ for each $v,v'\in|\Delta|$. We say that $h$
is strictly convex on $\Delta$ if moreover $h|\sigma\neq h|\sigma'$ for each $\sigma,\sigma'\in\Delta(l)$.
Observe that, given any $h\in SF(\Delta,M')$, there is a positive integer $n$ such that $nh$ is
$(\Delta,M)$-linear function. One can show that, if a torus $S'$ is a etale cover of $S$ with character group
$M'$, then we can identify $SF(\Delta,M')$ with the group of $S'$-linearized line bundles on $Z$ ($S'$ acts on
$Z$ by the quotient map $\xymatrix{ S'\ar@{->>}[r]&S  }$). We can suppose that $-h|\sigma$ is the character of
$S'$ on the fibre $L(x_{\sigma})$ over the $S$-stable point $x_{\sigma}$ associated to any maximal cone $\sigma$
of $\Delta$.

\begin{deff}
Let $L_{h}$ be a line bundle on $Y$. We will say that $L_{h}$ is
almost spherical  if $h|\sigma$ belongs to lattice generated by the
spherical weights for each $\sigma\in\Delta(l)$. Moreover we will
say that $L_{h}$ is spherical  if $h|\sigma$ is a spherical weight
for each $\sigma\in\Delta(l)$. We will say that $h$ is almost
spherical (respectively spherical) if $L_{h}$ is almost spherical
(respectively spherical).
\end{deff}

Observe that if $L_{h}$ is an almost spherical line bundle, then we
can think $h$ as a $(\Delta,\Lambda_{X})$ linear function and
$h^{c}$ as a $(\Delta^{c},\Lambda_{X})$ linear function. We now want
to describe the sections of a line bundle over $Y$. Observe that the
space of sections is multiplicity-free because $Y$ is a spherical
variety, i.e. it has a dense $B$-orbit. Now we want to define  sets
in bijective  correspondence with the sets of  the irreducible
subrepresentations respectively of $H^{0}(Y,L)$, $H^{0}(Z,L|Z)$ and
$H^{0}(Z^{c},L|Z^{c})$.

\begin{deff}
Given $h\in\Lambda_{Y}$ let
\[\Pi(Z,h)=\{\mu\in\bigcap_{\sigma\in\Delta(l)}(h|\sigma+(M\cap\check{\sigma}))\},\]
\[\Pi(Z^{c},h)=\{\mu\in\bigcap_{\sigma\in\Delta^{c}(l)}(h|\sigma+(M\cap\check{\sigma}))\]
and
\[\Pi(Y,h)= \Pi(Z,h)\cap \Lambda^{+}.\]
\end{deff}

Before we describe the sections of $L_{h}$, we want to rewrite the conditions for a weight to belong to
$\Pi(Z,h)$, respectively to $\Pi(Z^{c},h)$.

\begin{lem}\label{csv 2}
Let $\lambda$ be a weight in $\Lambda_{X}$ and let $h$ be in $\Lambda_{Y}$. Then:
\begin{enumerate}
\item $\lambda\in\Pi(Z,h)$ if and only if $\lambda\geq h$ as functions on
$|\Delta|$
\item $\lambda\in\Pi(Z^{c},h)$ if and only if $\lambda\geq h^{c}$ as functions on $N_{\textbf{R}}$
\end{enumerate}
\end{lem}

\begin{thm}[Theorem 3.4 in \cite{B}]\label{B 3.4}
Let $L_{h}$ be a line bundle on $Y$. Then
\[H^{0}(Y,L_{h})=\bigoplus_{\mu\in\Pi(Y,h)}V^{*}_{\mu}.\]
In particular $H^{0}(Y,L_{h})\neq0$ if and only if $\Pi(Y,h)$ is not empty.
\end{thm}

We want to give an idea of a possible construction of $H^{0}(Y,L_{h})$. We need the following lemma.
\begin{lem}[Lemma 2.7 in \cite{B}]\label{B 2.7} $O(Y_{\tau})$ is an almost spherical line bundle whose
associated $\Delta$-linear function  $d^{\tau}$ satisfies the equalities $d^{\tau}(\varrho(\tau'))=
-\delta_{\tau,\tau'}$ for each $\tau'\in\Delta(1)$. Moreover, if define $s_{\tau}\in H(Y,L_{d^{\tau}})$ as the
unique, up to a scalar,  section $s_{\tau}\in H(Y,L_{d^{\tau}})$ whose divisor is $Y_{\tau}$, then it is $G$
invariant.
\end{lem}
Observe that, given $\lambda\in\Lambda_{X}$, then $\lambda\in\Pi(Z,h)$ if and only if
$h=\lambda+\sum_{\tau\in\Delta(1)}a_{\tau}d^{\tau}$ where $a_{\tau}$ is a positive integer for each
$\tau\in\Delta(1)$. Let $\lambda\in\Pi(Y,h)$ and  remember that $H^{0}(X,L_{\lambda})$ contains
$V_{\lambda}^{*}$. $H^{0}(X,L_{\lambda})\subset H^{0}(Y,\pi^{*}L_{\lambda})$, so $H^{0}(Y,\pi^{*}L_{\lambda})$
contains a lowest weight vector $v_{-\lambda}$ of weight $-\lambda$. There are positive constants $a_{\tau}$
such that $h-\lambda=\sum_{\tau\in\Delta(1)}a_{\tau}d^{\tau}$ because of the lemma~\ref{csv 2}. Thus
$v_{-\lambda}\cdot\prod s_{\tau}^{a_{\tau}}$ is a not-zero  section of $H^{0}(Y,L_{h})$ with weight $-\lambda$
because the sections $s_{\tau}$ are $G$-invariant. Moreover $v_{-\lambda}\cdot\prod s_{\tau}^{a_{\tau}}$ is
invariant by  the unipotent part of the opposite $B^{-}$ of the fixed Borel group of $G$. Thus
$H^{0}(Y,L_{h})\supseteq\bigoplus_{\mu\in\Pi(Y,h)}V^{*}_{\mu}$. Because of the previous theorem we give the
following definition:

\begin{deff}
Given $h$  in $\Lambda_{Y}$ and  $\lambda$  in $\Pi(Y,h)$, we write
$h=\lambda+\sum_{\tau\in\Delta(1)}a_{\tau}d^{\tau}$ for suitable $a_{\tau}\in \textbf{Z}^{+}$. We define
$s^{h-\lambda}$ as the section $\prod s_{\tau}^{a_{\tau}}$ of $H^{0}(Y,L_{h-\lambda})$.
\end{deff}

We want to describe also the the sections over $Z$, respectively over $Z^{c}$.

\begin{prop}
\label{B 4.1} Let $L_{h}$ be a line bundle on $Y$. Then
\begin{enumerate}
\item \[H^{0}(Z,L_{h}|Z)=\bigoplus_{\mu\in\Pi(Z,h)}\textbf{C}\chi^{\mu},\]
where $\chi^{\mu}$ is a $T$-seminvariant section of weight $-\mu$. In particular $H^{0}(Z,L_{h}|Z)$ $\neq0$ if
and only if $\Pi(Z,h)$ is not empty.

\item \[H^{0}(Z^{c},L_{h}|Z^{c})=\bigoplus_{\mu\in\Pi(Z^{c},h)}\textbf{C}\chi^{\mu}.\]
In particular $H^{0}(Z^{c},L_{h}|Z^{c})\neq0$ if and only if $\Pi(Z^{c},h)$ is not empty.

\end{enumerate}
\end{prop}

{\em Remark.} Let $\pi:Y\rightarrow Y'$ be an $G$-equivariant morphism between two complete symmetric varieties
and let $L_{h}$ be a line bundle on $Y'$. Then the pullback $\pi^{*}(L_{h})$ is the line bundle on $Y$
associated to $h$ and $H^{0}(Y,\pi^{*}(L_{h}))= H^{0}(Y',L_{h})$.

Now we want to explain some relations between the previous sets.

\begin{cor}[Corollary 4.1 in \cite{B}]\label{B 4.1 cor}
Given   $h\in \Gamma_{Y}$, we have the equality

\[\Pi(Y,h)= \Pi(Z^{c},h)\cap \Lambda^{+} .\]
\end{cor}

If $L_{h}$ is almost spherical, we can say more.

\begin{prop}[Proposition 4.2 and theorem 4.2 in \cite{B}]
\label{B 4.2} If $h\in \Lambda_{Y}$ is almost spherical, then

\[\Pi(Z^{c},h)=\bigcup_{w\in W^{1}} w\cdot\Pi(Y,h).\]
Moreover the restriction map $H^{0}(Y,L_{h})\rightarrow H^{0}(Z^{c},L_{h}|Z^{c})$ is surjective.
\end{prop}

We want to make some remark on the second part of the proposition. Let $w\cdot\mu\in\Pi(Z^{c},h)$ with
$\mu\in\Pi(Y,h)$. Any section $s\in V_{\mu}^{*}\subset H^{0}(Y,L_{h})$  of  weight $-w\cdot\mu$ has not-zero
restriction to $Z^{c}$ because $U^{-}\cdot Z^{c}$ is dense in $Y$. Moreover, up to choose another basis of the
root system, we can suppose that it is lowest weight vector.

Remember  that there is an one-to-one  correspondence between the convex functions on $N_{\textbf{R}}$ with
values in $\textbf{R}\cup\infty$ and the convex sets in $M_{\textbf{R}}$, which send a convex function $h$ to
the convex set $Q_{h}=\{m\in M_{\textbf{R}}: m(v)\geq h(v),\ \forall\ v\in N_{\textbf{R}}\}$. Moreover
$h(v)=inf\{m(v); \ m\in Q_{h}\}$  for each $ v\in N_{\textbf{R}}$. To every linearized line bundle $L_{h}$ on a
toric variety we associated the polyhedron $Q_{h}$. We want to do the same with the line bundles on a complete
symmetric variety. (Remember that a polyhedron is the intersection of a finite number of semispaces).

\begin{deff}
Let $Y$ be a complete symmetric variety and let $L_{h}$ be an almost spherical line bundle on $Y$.  We define
the polytope associated to $L_{h}$ as

\[P_{h}=\{m\in M_{\textbf{R}}: m(v)\geq h^{c}(v)\ \forall  v\in |\Delta^{c}| \}.\]

Moreover we define the polyhedron associated to $h$ as

\[Q_{h}=\{m\in M_{\textbf{R}}: m(v)\geq h(v)\ \forall v\in |\Delta| \}.\]
\end{deff}
Observe that $P_{h}=\{m\in M_{\textbf{R}}: m(\varrho(\tau))\geq h^{c}(\varrho(\tau))\ \forall \ \tau\in
\Delta^{c}(1) \}$ and $Q_{h}=\{m\in M_{\textbf{R}}: m(\varrho(\tau))\geq h(\varrho(\tau))\ \forall \tau\in
\Delta(1) \}$.

\section{Ample line bundles and line bundles generated by global sections}

Brion \cite{Br} has found a characterization of the ample line bundles (respectively the line bundles generated
by global sections) on a spherical variety. Now we want to find different   conditions for a line bundle on a
complete symmetric variety to be generated by global sections, respectively to be ample.

\begin{prop}
\label{ampio} Let $L_{h}$ be a line bundle on $Y$. Then

\begin{enumerate}
\item $L_{h}$ is generated by global sections if and only if $h$ is convex and
$h|\sigma$ is dominant for each $\sigma\in\Delta(l)$.

\item $L_{h}$ is very ample if and only if $h$ is strictly convex on
$\Delta$ and  $h|\sigma$ is a regular  weight for each $\sigma\in\Delta(l)$.

\item $L_{h}$ is ample if and only if it is very ample.
\end{enumerate}
\end{prop}

{\em Proof.} One can easily show the necessity of the conditions.
Indeed if $L$ is generated by global sections (respectively is
ample) then the restrictions of $L$ to $Z$ and to any closed orbit
$O_{\sigma}$ are generated by global sections (respectively are
ample). To prove the sufficiency of the condition in the first point
we need a lemma.

\begin{lem}\label{am l1}
If $h$ is convex and $h|\sigma$ is dominant then the restriction map to the closed orbit $O_{\sigma}$
\[H^{0}(Y,L_{h})\rightarrow H^{0}(O_{\sigma},L_{h}|O_{\sigma})\]
is surjective.
\end{lem}

{\em Proof.} We have
$H^{0}(O_{\sigma},L_{h}|O_{\sigma})=V_{h|\sigma}^{*}$ because
$h|\sigma$ is dominant. Moreover there is a lowest weight vector
$\varphi\in  H^{0}(Y,L_{h})$ of weight $-h|\sigma$ because of the
convexity of $h$. Hence, because of the irreducibility of
$V_{h|\sigma}^{*}$, it is sufficient to prove that the restriction
of $\varphi$ to $O_{\sigma}$ is not zero. Observe that $\varphi=
\varphi'\cdot\prod s_{\tau_{i}}^{a_{i}}$, where $\varphi'$ is a
lowest weight vector of $ V_{\lambda_{\sigma}}^{*}\subset
H^{0}(Y,L_{h|\sigma})$ and $a_{i}>0$ only if $\tau_{i}$ is not
contained in $\sigma$. Recall the $s_{\tau_{i}}$ vanishes exactly on
$Y_{\tau_{i}}$ and observe that $\varphi'|O_{\sigma}\neq0$ because
$L_{h|\sigma}|O_{\sigma}$ is the line bundle corresponding to
$h|\sigma$. $\square$

Now we can prove the sufficiency of the condition in the first point. If the locus of  base points is not empty,
then it contains a closed orbit $O_{\sigma}$, because it is closed and stabilized by $G$. Given any $y\in
O_{\sigma}$ there is a section $\widetilde{s}\in H^{0}(O_{\sigma},L_{h}|O_{\sigma})$ such that
$\widetilde{s}(y)\neq0$ because $h|\sigma$ is dominant. This section can be extended to a global section over
$Y$  because of the previous lemma. This is a contradiction.

Now we want to show the sufficiency of the condition in the second point if $L_{h}$ is  spherical. First of all
we want to show that $L_{h}|Z$ is very ample, or better still that  $L_{h}|Z^{c}$ is very ample. $h^{c}$ is
convex because of the first point of the proposition.  Suppose by contradiction that there are two distinct
maximal cones $\sigma$ and $\sigma'$ such that $h|\sigma= h|\sigma'$.  Let $\gamma$ be the convex cone formed by
the points $v$ such that $h(v)=(h|\sigma)(v)$. We can suppose that $\sigma\in\Delta$ and that
$\sigma'\nin\,\Delta$, so there is an hyperplane $H$  secant to $\gamma$ and whose intersection with
$\sigma(e_{1},...,e_{l})$ is a face of $\sigma(e_{1},...,e_{l})$, say
$\sigma(e_{1},...,\widehat{e}_{i},...,e_{l})$. Moreover there is a vector $v$ of $\gamma$ contained in the
interior of $|\Delta|$ and such that also $s_{i}v$ belongs to $\gamma$.

By hypothesis we know that $h^{c}(v) =(h|\sigma)(v)$ and that $h^{c}(s_{i}v)= (s_{i}\cdot(h|\sigma))(v)$, so
$(h|\sigma)(v)=(s_{i}(h|\sigma))(v)$ because of the invariance of $h^{c}$ by $W^{1}$.
$((h|\sigma)-s_{i}(h|\sigma))(v)$ is strictly positive because $h|\sigma$ is a regular weight and $v$ is in the
interior of $|\Delta|$; this is a contradiction. Observe that we have proved a more general statement. Let $Z$
be  a possibly singular toric variety proper over $\textbf{A}^{l}$ and let $L_{h^{c}}$ be a linearized line
bundle on $Z^{c}$ such that $h^{c}$ is invariant for the action of $W^{1}$. If, for each $\sigma\in\Delta(l)$,
$h^{c}|\sigma$ is a regular weight (respectively a dominant weight) and $h$ is strictly convex (respectively
convex) on the fan of $Z$, then $L_{h^{c}}$ is ample (respectively generated by global sections).

Since $L_{h}$ is generated by global sections, we have an
equivariant morphism $\varphi: Y \rightarrow \textbf{P}(V)$ with
$V=H^{0}(Y,L)^{*}$. Let  $U$ be the locus where $\varphi$ is not an
embedding. We could try to prove this point like the previous one,
namely using the fact that the restriction of $L_{h}$ to $Z$,
respectively to any closed orbit is very ample. Instead we will use
the stronger fact that $L_{h}|Z^{c}$ is ample. Observe that the
restriction of the sections from $Y$ to $Z^{c}$ is surjective, while
the restriction of the sections to $Z$  is clearly not surjective.

Now we want to show that $U$ is stable and closed in the Euclidean
topology. $U$ is the union of two loci: the locus  $U_{1}$ of the
points where the differential of $\varphi$ is not injective and the
locus $U_{2}$ of the points where $\varphi$ is not injective.
$U_{1}$ and $U_{2}$ are $G$-stable because $\varphi$ is equivariant.
$U_{1}$ is closed because it is the locus of the zeroes of the
jacobian of $\varphi$. Now we want to prove that the closure of
$U_{2}$ is contained in $U$. Let $\{x_{n}\}$ be any sequence in
$U_{2}$ and suppose that it converges to $x\in Y$.   By hypothesis
there is a sequence $\{y_{n}\}$ in $U_{2}$ such that $x_{n}\neq
y_{n}$ and  $\varphi(x_{n})=\varphi(y_{n})$ for each $n$. We can
suppose,  up to taking sub-sequences, that $\{y_{n}\}$ has limit $y$
in $Y$ and that $\varphi(x)=  \varphi(y)$. Suppose by contradiction
that $x$ does not belong to $U$, in particular $x=y$. Because of the
Dini theorem there is a open neighborhood $W$ of $x$ such that
$\varphi|W$ is a diffeomorphism onto the image $\varphi(W)$. This is
a contradiction because there is an integer $n_{0}$ such that
$x_{n}$ and $y_{n}$ belong to $W$ for each $n>n_{0}$. In particular,
we have proved that if $U_{1}$ is empty then $x$ must be different
from $y$.

First suppose that $U_{1}$ is not empty, so it contains a closed orbit $O_{\sigma}$. Let  $x_{\sigma}$ be the
intersection of $Z$ and $O_{\sigma}$. The map $\xymatrix{H^{0}(Z^{c},L_{h}|Z^{c})^{*}\ar[r]&H^{0}(Y,L_{h})^{*}}$
dual to the restriction map is injective, so we have a commutative diagram

\[\xymatrix{Y \ar[rr]^{\varphi}& &\textbf{P}(H^{0}(Y,L_{h})^{*})\\Z^{c}\ar[u]
\ar[rr]^{\varphi'}&& \textbf{P}(H^{0}(Z^{c},L_{h}|Z^{c})^{*})\ .\ar[u]}\] $\varphi'$ is an embedding because
$h^{c}$ is strictly convex on $\Delta^{c}$, so $\varphi(Z^{c})$ is isomorphic to $Z^{c}$. Let $[h]$ be the image
$\varphi(H)$ of $H\in G/H$ and let $[v_{h|\sigma}]$ be the image of $x_{\sigma}$. Observe that $[v_{h|\sigma}]$
is the class of a highest weight vector of $V_{h|\sigma}\subset
H^{0}(Y,L_{h})^{*}$. 
We can write $H^{0}(Y,L_{h})^{*}=V_{h|\sigma}\oplus V'$ for a suitable representation $V'$.  We can choose
$h=v_{h|\sigma}+\sum v_{i}$ where the $v_{i}$ are weight vectors with weights contained in $h|\sigma+M$.  Let
$\widetilde{A}$ be the affine open set of $\textbf{P}(H^{0}(Y,L_{h} )^{*})$ where a lowest weight vector $s\in
H^{0}(Y,L_{h})$ of weight $-h|\sigma$ is not zero, i.e. $\widetilde{A}=v_{h|\sigma}+V'_{h|\sigma}\oplus V'$
where $V_{h|\sigma}=\textbf{C}v_{h|\sigma}\oplus V'_{h|\sigma}$ in a $T$-equivariant way. $A:=\widetilde{A}\cap
\varphi(Y)$  is $B^{-}$ stable and its intersection with $\varphi(Z^{c})$ is $\varphi(U_{\sigma})$, so
$\varphi(x_{\sigma})$ belongs to $A$. Indeed the set of points $\{x\in Z^{c}:\ s(x)=0\}$ is the union of the
divisor $Z_{\tau}^{c}$ for $\tau\nsubseteq\sigma$.

We want to study the restriction of $\varphi$ to $U^{-}\cdot U_{\sigma}$, where $U^{-}$ is the unipotent group
whose Lie algebra is $\bigoplus_{\alpha\in -\phi_{1}^{+}}\mathfrak{g}_{\alpha}$. $H^{0}(Y,L_{h})^{*}$ is
isomorphic to its dual in a $\theta$ linear way (see lemma 1.6~\cite{CSV I}). Thus there is a not degenerate
bilinear form $(\, , )$ on $H^{0}(Y,L_{h})^{*}$ with the following properties: 1) given any two distinct
irreducible components $V_{1}$ and $V_{2}$ of $H^{0}(Y,L_{h})^{*}$, they are orthogonal; 2)
$(gu,v)=(u,\theta(g^{-1})v)$ for each $g\in G$, $u,v \in H^{0}(Y,L_{h})^{*}$. Let $\Upsilon'$ be the tangent
space in $v_{h|\sigma}$ to the orbit $U^{-}\cdot v_{h|\sigma}$ and let $\Upsilon$ be the space generated by
$\Upsilon'$ and $v_{h|\sigma}$. One can show that the restriction of $(\ , \ )$ to $\Upsilon$ is non-degenerate,
that $\Upsilon$ is stable under $P$ and that the orthogonal $\Upsilon^{\bot}$ is stable under $\theta(P)$ (see
lemma 2.4~\cite{CSV I}). Moreover $U_{\sigma}\subset v_{h|\sigma}+\Upsilon^{\bot}$ (see lemma 2.5 in~\cite{CSV
I}). Let $q$ be the projection of $H^{0}(Z^{c},L_{h}|Z^{c})^{*}$ onto
$H^{0}(Z^{c},L_{h}|Z^{c})^{*}/\Upsilon^{\bot}$. $U^{-}\subset \theta(P)$, so $U^{-}$ acts on
$H^{0}(Z^{c},L_{h}|Z^{c})^{*}/\Upsilon^{\bot}$ and the projection is equivariant.   The affine hyperplane
$q(\widetilde{A})$ in $H^{0}(Z^{c},L_{h}|Z^{c})^{*}/\Upsilon^{\bot}$ is stable by the action of $U^{-}$. One can
show that the map $j:U^{-}\rightarrow\pi(\widetilde{A})$ defined by $j(u)=q(uv_{h|\sigma})$ is an $U^{-}$
equivariant isomorphism (see lemma 2.6 in~\cite{CSV I}). Thus the tangent space to $U_{\sigma}$ at
$v_{h|\sigma}$ is orthogonal to $\Upsilon$ and the differential of $\varphi$ is injective in $x_{\sigma}$. Hence
$U_{1}=\emptyset$, so $U_{2}$ is equal to $U$ and it is closed.

Now suppose that $U_{2}$ is not empty, so it contains a closed orbit $O_{\sigma}$. Given $x\in O_{\sigma}$,
there is $y\neq x$ such that $\varphi(x)=\varphi(y)$. We can suppose that $y$ belongs to a closed orbit. Indeed
there is an element $g$ of $G$ and an one parameter subgroup $\gamma$ of $T$ such that $y'=lim_{t\rightarrow
0}\gamma(t)gy$ is a point of $Z$ fixed by $T$, in particular $y'$ belongs to a closed orbit $O_{\sigma'}$.
Moreover $x'=lim_{t\rightarrow 0}\gamma(t)gx$ belongs to $O_{\sigma}$ and $\varphi(x')=\varphi(y')$. By the
previous part of the proof $x'$ is different by $y'$.

The closed orbits $O_{\sigma}$ and $O_{\sigma'}$ are different because  $L|O_{\sigma}$ is very ample. Because of
the  lemma~\ref{am l1} there is a global section $s$, lowest weight vector of weight $-h|\sigma$, which does not
vanish on $O_{\sigma}$, so   we can suppose that $s(x)\neq0$ up to a translation. Since $s$ is $U^{-}$ invariant
and $h$ is strictly convex on $\Delta$, $s$ vanishes on a $G$-stable divisor of $Y$ which contains
$O_{\sigma'}$. Hence we have obtained a contradiction, namely $\varphi(x)\neq\varphi(y)$.

Finally we can consider the exceptional case. First of all we want to recall some facts. Let $Y$ be a complete
exceptional symmetric variety and let $X$ be the corresponding wonderful variety. We can chose an order of the
simple roots of $\phi$ such that $\alpha_{1},...,\alpha_{s}$ are exceptional roots with the following property:
$Pic(X)$ is generated by the spherical weights and by the fundamental weights
$\omega_{\alpha_{1}},...,\omega_{\alpha_{s}}$ corresponding respectively to $\alpha_{1},...,\alpha_{s}$.
Moreover,  given an element $h$ of $\Lambda_{Y}$ such that $h|\sigma$ is dominant for each $\sigma\in\Delta(l)$,
there  are  integers $a_{i}$ such that $h'=h-\sum a_{i}\omega_{\alpha_{i}}$ is  another element  of
$\Lambda_{Y}$ and $h'|\sigma\in\Omega$ for each $\sigma\in\Delta(l)$. For each $i$ we can suppose that $a_{i}$
is positive, up to exchange $\alpha_{i}$ with $\overline{\theta}(\alpha_{i})$. If $L_{h}$ satisfies the
hypotheses of the second point, then $L_{h'}$ is very ample and $L_{h-h'}$ is generated by global sections, so
$L_{h}$ is very ample.

The third point is obvious. $\square$

Now, we can characterize the ample  line bundles  on $Z$.

\begin{prop}\label{tor lb5}
Suppose that $Z$ is a (possibly singular) toric variety proper over $\textbf{A}^{l}$ and let $h\in
SF(\Delta,M)$. Then $L_{h}$ is ample  if and only if   $h$ is strictly convex on $\Delta$. If $Z$ is smooth then
$L_{h}$ is ample  if and only if it is very ample.
\end{prop}

{\em Proof.} Suppose that $L_{h}$ is ample, then there is an integer $n$ such that $L_{nh}$ is very ample, in
particular   $nh$ is convex. Hence $L_{nh}$ is the pullback of a line bundle generated by global sections on a
variety $Z'$ such that  $Z$ is proper over $Z'$ and $nh$ is strictly convex on the fan $\Delta'$ of $Z'$. Let
$\varphi:Z\rightarrow \textbf{P}(V)$ be an embedding such that $L_{nh}=\varphi^{*}\mathcal{O}(1)$, then
$\varphi$ factorizes through $Z'$ because $H^{0}(Z',L_{h})=H^{0}(Z,L_{h})$. Since $\varphi$ is an embedding,
$Z'$ must be $Z$, so $h$ is strictly convex on $\Delta$. The viceversa is implied by the following lemma.

\begin{lem}\label{csv 6a}
Let $Z$  be is a (possibly singular) toric variety proper over
$\textbf{A}^{l}$ and let $L$ be a line bundle on $Z$ generated by
global sections.  Given any homogeneous symmetric variety $G/H$ of
rank $l$, let $Z^{c}$ be the complete toric variety associated to
$Z$ and observe that $N_{H^{0}}(S)$ acts on $Z^{c}$. Then there is a
linearized line bundle $L'$ on $Z^{c}$ generated by global sections
and a  $N_{H^{0}}(S)$ linearization of $L'$  such that: i) the
restriction of $L'$ to $Z$ is $L$; 2) for each $S$-fixpoint $x$ of
$Z$,  the character of $S$ on the fibre $L'(x)$ is a regular weight.
\end{lem}

{\em Proof.} We can suppose that $Z$ is smooth up to take  a resolution of singularities  of $Z$. Given any
homogeneous symmetric variety $G/H$ of rank $l$, let $Y$ be the complete symmetric variety   associated to $Z$.
Recall that there is almost spherical line bundle $L_{h}$ on $Y$ whose restriction to $Z$ is $L$ and let
$\lambda$ be a regular spherical weight. There is a positive integer $n$ such that $(h+n\lambda)|\sigma$ is a
regular weight for each $\sigma\in\Delta(l)$, so $L'=L_{h+n\lambda}$ satisfies ours requests. $\square$.

We now  can conclude the proof the proposition~\ref{tor lb5}.
Observe that, up to changing the linearization of $L_{h}$, we can
suppose that the line bundle $L_{h^{c}}$ is as in the previous
lemma. Moreover $L_{h^{c}}$ is ample on $Z^{c}$ and $L_{h}$ is ample
on $Z$ because of the proof of the proposition~\ref{ampio}. The last
point of the proposition is implied by the Demazure theorem (see
\cite{De}). $\square$

\begin{cor}
Let $Y$ be a complete symmetric variety. The following conditions are equivalent:
\begin{enumerate}
\item $Y$ is projective;

\item $Z^{c}$ is projective;

\item $Z$ is quasiprojective.
\end{enumerate}
\end{cor}

Now we want to reformulate the proposition~\ref{ampio} using $h^{c}$ instead of $h$.

\begin{prop}\label{ampio 2} Let  $h$ be an almost spherical $\Delta$-linear function, then

\begin{enumerate}
\item $h^{c}$ is convex on $\Delta^{c}$ if and only if $h$ is convex on $\Delta$ and $h|\sigma$ is
dominant for each $\sigma\in\Delta(l)$.

\item $h^{c}$ is strictly convex on  $\Delta^{c}$ if and only if
$h$ is strictly convex on  $\Delta$ and $h|\sigma$ is a regular weight for each $\sigma\in\Delta(l)$.
\end{enumerate}
\end{prop}

{\em Proof.} By the proof of the proposition~\ref{ampio} it is sufficient to prove the following facts. If
$h^{c}$ is an almost spherical convex $\Delta^{c}$ linear function then $h$ is a spherical $\Delta$-linear
function. If $h^{c}$ is also strictly convex on $\Delta^{c}$, then $h|\sigma$ is regular for each
$\sigma\in\Delta(l)$.

Given $\sigma\in\Delta(l)$, there is an element $w\in W^{1}$ such that $w\cdot h|\sigma$ is a dominant weight.
$h|\sigma-w\cdot h|\sigma$ has  positive values on $|\Delta|$ and $w\cdot h|\sigma$ is the restriction of
$h^{c}$ to $w^{-1}\cdot\sigma$. Let  $v$ be a vector in the interior of $\sigma$, then $(w\cdot
h|\sigma)(v)=(h|\sigma)(v)$ because of the convexity of $h^{c}$. We have $w\cdot h|\sigma=h|\sigma$ because $v$
is a vector inside the Weyl chamber $|\Delta|$, so $h|\sigma$ is dominant. If $h^{c}$ is strictly dominant on
$\Delta^{c}$, then $h|\sigma$ is different from $w\cdot h|\sigma$ for each $w\in W^{1}$, so $h|\sigma$ is
regular.  $\square$

\section{Reduction to the complete toric variety}

\textit{In the following we will always suppose that $L_{h}$ is an almost spherical line bundle, unless we will
explicitly say otherwise.} We start to study the multiplication of sections of two line bundles on $Y$. First of
all, we want to show that this problem is equivalent to the similar problem on the complete toric variety
$Z^{c}$ associated to $Y$. Let $L_{h}$ and $L_{k}$ be any two line bundles on $Y$ generated by global sections.
Let

\[ M_{h,k}:H^{0}(Y,L_{h})\otimes H^{0}(Y,L_{k})\longrightarrow
H^{0}(Y,L_{h+k})\] be the product of sections on $Y$ and let
\[ m^{c}_{h,k}:H^{0}(Z^{c},L_{h}|Z^{c})\otimes H^{0}(Z^{c},L_{k}|Z^{c})\longrightarrow
H^{0}(Z^{c},L_{h+k}|Z^{c})\] be the product of sections of the restrictions to $Z^{c}$ of these line bundles.

\begin{thm}\label{red ctor}
Let $L_{h}$ and $L_{k}$ be two almost spherical line bundle on $Y$ generated by global sections. Then $M_{h,k}$
is surjective if and only if $m^{c}_{h,k}$ is surjective.
\end{thm}

{\em Proof.} The necessity of the condition is implied by the surjectivity of the restriction maps from $Y$ to
$Z^{c}$. Indeed if $i:Z^{c}\rightarrow Y$ is the canonic inclusion then $m^{c}_{h,k}\circ (i^{*}\otimes
i^{*})=i^{*}\circ M_{h,k}$. Now suppose that $m^{c}_{h,k}$ is surjective. It is sufficient to show that the
image of $M_{h,k}$ contains a basis of semi-invariant sections. If $h$ and $k$ are linear then  $M_{h,k}$ is
surjective by the theorem~\ref{CM}. In general, given $\nu\in \Pi(Y,h+k)$ there are $\lambda\in\Pi(Z^{c},h)$ and
$\mu\in\Pi(Z^{c},k)$ such that $\nu=\lambda+\mu$. Moreover there are elements $w_{1}$ and $w_{2}$ of $W^{1}$
such that $w_{1}\cdot\lambda$ and $w_{2}\cdot\mu$ are dominant weights. Observe that $\nu\geq w_{1}\cdot\lambda+
w_{2}\cdot\mu$ on $|\Delta|$. Moreover $w_{1}\cdot\lambda\geq h$ and $w_{2}\cdot\mu\geq k$ because $h^{c}$ and
$k^{c}$ are convex and invariant for the action of $W^{1}$. Thus
$s^{h-w_{1}\cdot\lambda}H^{0}(Y,L_{w_{1}\cdot\lambda})\subset H^{0}(Y,L_{h})$ and
$s^{k-w_{2}\cdot\mu}H^{0}(Y,L_{w_{2}\cdot\mu})\subset H^{0}(Y,L_{k})$. We know that $Im\,M_{w_{1}\cdot\lambda,
w_{2}\cdot\mu}$ contains a lowest weight vector  $\varphi\in H^{0}(Y,L_{w_{1}\cdot\lambda+ w_{2}\cdot\mu})$  of
weight $-\nu$.  Thus $s^{h+k-w_{1}\cdot\lambda- w_{2}\cdot\mu}\varphi$ is contained in
$s^{h+k-w_{1}\cdot\lambda- w_{2}\cdot\mu}Im\,M_{w_{1}\cdot\lambda, w_{2}\cdot\mu}\subset Im\,M_{h,k}$ and it is
not zero. $\square$

We can prove the following proposition  without assuming  the surjectivity of $m^{c}_{h,k}$. Given   two convex
$(\Lambda_{X},\Delta)$ linear function, say $h$ and $k$, let  $\Pi(Y,h,k)$ be the set of the weights of the
lowest weight vectors contained in $Im M_{h,k}$.

\begin{prop}\label{red ctor2}
$\Pi(Y,h,k)$ is saturated with respect to the dominant order of the roots in  $\widetilde{\phi}$.
\end{prop}

{\em Proof.} $\Pi(Y,h+k)$ is saturated because the simple restricted roots have negative values on $|\Delta|$.
Given $\nu\in \Pi(Y,h,k)$ there are two  weights  $\lambda\in\Pi(Y,h)$ and $\mu\in\Pi(Y,k)$ such that
$\nu=\lambda+\mu$.  Moreover there are element $w_{1}$, $w_{2}$ in the Weyl group $W^{1}$ such that
$w_{1}\cdot\lambda$ and $w_{2}\cdot\mu$ are dominant weights. Observe that $\nu\geq w_{1}\cdot\lambda+
w_{2}\cdot\mu$ on $|\Delta|$, so $\nu\in\prod(Y,w_{1}\cdot\lambda+w_{2}\cdot\mu)$. Let   $\nu'$ be a spherical
weight dominated by $\nu$, then $\nu'\in \Pi(Y,w_{1}\cdot\lambda+w_{2}\cdot\mu)$ because this set is saturated.
Let  $\varphi$ be a lowest weight vector of weight $\nu'$. Because of the surjectivity of
$M_{w_{1}\cdot\lambda,w_{2}\cdot\mu}$ we have $\varphi\in s^{h+k-w_{1}\cdot\lambda-w_{2}\cdot\mu}Im
M_{w_{1}\cdot\lambda,w_{2}\cdot\mu}\subset Im M_{h,k}$. $\square$

\section{Reduction to the open toric variety}

In this section  we want to show that, given two ample line bundles on $Y$, the product of sections on $Z^{c}$
is surjective if and only if the product of sections on $Z$ is surjective. Moreover we will study the relation
between the sections of $L|Z$ and the sections of $L|Z^{c}$ for any ample line bundle $L$ on $Y$. Before we have
to define  some notations. We  fix a cone $\sigma\in\Delta(l)$ and we  set $v_{h}=h|\sigma$ for each $h\in
SF(\Delta,\Lambda_{X})$, so $\prod(Z,h)=Q_{h}\cap(M+v_{h})$ and $\Pi(Z^{c},h)=P_{h}\cap(M+v_{h})$. Given any
one-dimensional cone $\tau$, we set $\rho(\tau)$ as the primitive vector of $\tau$ i.e. the  more little
not-zero vector of $N\cap \tau$.

$\{e_{1},...,e_{l}\}$ is the basis of $N_{\textbf{R}}$ dual to the
basis $\{f_{1},...,f_{l}\}$ of $M_{\textbf{R}}$. We have to define a
second basis $\{g_{1},...,g_{l}\}$ of $M_{\textbf{R}}$ because the
fundamental Weyl chamber $C^{+}$ is more easily defined using the
basis the fundamental weights than the basis of the simple roots.
$g_{i}$ is a positive multiple of   $-\omega_{i}$, more precisely
$-g_{i}$ is the $i$-th fundamental weight of the unique reduced root
system contained in $\widetilde{\phi}$ which share a basis with
$\widetilde{\phi}$. $g_{1},...,g_{l}$ generate a lattice which
contains $M$.  Let $\{\check{g}_{1},...,\check{g}_{l}\}$ be the dual
basis of $\{g_{1},...,g_{l}\}$. Observe that
$C^{+}=\sigma(-g_{1},...,-g_{l})$. Given a point $p$ in
$M_{\textbf{R}}$ we will use the following notations: $p=\sum
x_{i}f_{i}=\sum \dot{x}_{i}g_{i}$, using the "normal" coordinates
for the basis $\{f_{1},...,f_{l}\}$ and the "dotted" coordinates for
the basis $\{g_{1},...,g_{l}\}$. (In the following figures we
consider the case in which the restricted root system is of type
$A_{2}$ and $Z$ is $\textbf{A}^{2}$).

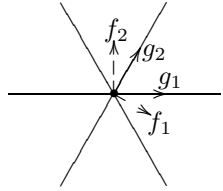
\[\begin{xy} <4pt, 0pt>:(0,0)*{\scriptstyle\bullet}="o";
"o"+0;"o"+(5,0)**\dir{-}?>*\dir{>}; "o"+0;"o"+(10,0)**\dir{-}; "o"+0;"o"+(-10,0)**\dir{-};
"o"+0;"o"+(2.5,4.3330)**\dir{-}?>*\dir{>}; "o"+0;"o"+(5,8.66)**\dir{-}; "o"+0;"o"+(-5,-8.66)**\dir{-};
 "o"+0;"o"+(-5,8.66)**\dir{-}; "o"+0;"o"+(5,-8.66)**\dir{-};
"o"+0;"o"+(3.4640,-2.0)**\dir{--}?>*\dir{>};"o"+0;"o"+(0,5)**\dir{--}?>*\dir{>};
"o"+(5.3,1.2)*{g_{1}};"o"+(3.8,3.85)*{g_{2}}; "o"+(0.2,6)*{f_{2}};"o"+(4.3,-2.77312)*{f_{1}};
\end{xy}\]

\begin{prop}\label{red otor l1a}
Let  $L_{h}$ be an ample spherical line bundle on $Y$. Then $Q_{h}\cap C^{+}=P_{h}\cap C^{+}$ and
$Q_{h}=P_{h}\cap C^{+}+\sigma(f_{1},...,f_{l})$.
\end{prop}

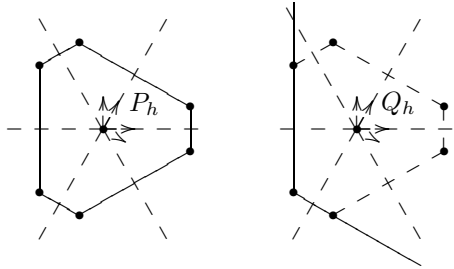
\[\begin{xy} <2.4pt, 0pt>:(0,0)*{\scriptstyle\bullet}="o";
(-10,-10)*{\scriptstyle\bullet}="p1"; (-10,10)*{\scriptstyle\bullet}="p2";
(-3.7648,13.6)*{\scriptstyle\bullet}="p3";"p3"+(17.4932,-10.1)*{\scriptstyle\bullet}="p4";
 (-3.7648,-13.6)*{\scriptstyle\bullet}="p8"; "p8"+(17.4932,10.1)*{\scriptstyle\bullet}="p7";
"o"+0;"o"+(5,0)**\dir{--}?>*\dir{>}; "o"+0;"o"+(15,0)**\dir{--}; "o"+0;"o"+(-15,0)**\dir{--};
"o"+0;"o"+(2.5,4.3330)**\dir{--}?>*\dir{>}; "o"+0;"o"+(10,17.320)**\dir{--}; "o"+0;"o"+(-10,-17.320)**\dir{--};
 "o"+0;"o"+(-10,17.320)**\dir{--}; "o"+0;"o"+(10,-17.320)**\dir{--};
"o"+0;"o"+(3.4640,-2.0)**\dir{--}?>*\dir{>};"o"+0;"o"+(0,5)**\dir{--}?>*\dir{>};
 "p1"+0;"p2"+0**\dir{-}; "p2"+0;"p3"+0**\dir{-}; "p3"+0;"p4"+0**\dir{-};
  "p4"+0;"p7"+0**\dir{-}; "p7"+0;"p8"+0**\dir{-}; "p8"+0;"p1"+0**\dir{-};
"o"+(6.5,4)*{P_{h}};
(40,0)*{\scriptstyle\bullet}="bo"; (30,-10)*{\scriptstyle\bullet}="bp1"; (30,10)*{\scriptstyle\bullet}="bp2";
(36.2352,13.6)*{\scriptstyle\bullet}="bp3";"bp3"+(17.4932,-10.1)*{\scriptstyle\bullet}="bp4";
(36.2352,-13.6)*{\scriptstyle\bullet}="bp8";"bp8"+(17.4932,10.1)*{\scriptstyle\bullet}="bp7";
"bo"+0;"bo"+(5,0)**\dir{--}?>*\dir{>}; "bo"+0;"bo"+(15,0)**\dir{--}; "bo"+0;"bo"+(-12,0)**\dir{--};
"bo"+0;"bo"+(2.5,4.3330)**\dir{--}?>*\dir{>}; "bo"+0;"bo"+(10,17.320)**\dir{--};
"bo"+0;"bo"+(-10,-17.320)**\dir{--};
 "bo"+0;"bo"+(-11,19.052)**\dir{--}; "bo"+0;"bo"+(10,-17.320)**\dir{--};
"bo"+0;"bo"+(3.4640,-2.0)**\dir{--}?>*\dir{>};"bo"+0;"bo"+(0,5)**\dir{--}?>*\dir{>};
 "bp1"+0;"bp2"+(0,10)**\dir{-}; "bp2"+0;"bp3"+0**\dir{--}; "bp3"+0;"bp4"+0**\dir{--}; "bp4"+0;"bp7"+0**\dir{--};
   "bp7"+0;"bp8"+0**\dir{--}; "bp1"+0;"bp8"+(13.856,-8)**\dir{-};
"bo"+(6.5,4)*{Q_{h}};
\end{xy}\]

The equations of $Q_{h}$ are of the form $\sum b_{i}x_{i}\geq b$ where the $b_{i}$ are positive constants. Thus
$Q_{h}$ is stable by translation with respect to vectors in $\sigma(f_{1},...,f_{l})$, i.e. $Q_{h}+\bigoplus
\textbf{R}^{+}(f_{i})\subset Q_{h}$. Let $H_{j}$ be the hyperplane of $M_{\textbf{R}}$ generated by
$g_{1},...,\widehat{g}_{j},...,g_{l}$, so the intersection of $H_{j}$ and $C^{+}$ is a Weyl wall. Let $s_{j}$ be
the orthogonal reflection with respect to $H_{j}$. Observe that, if $P_{h}$ contains a point $p$, then it
contains all the translates of $p$ by $W^{1}$, so it contains the orthogonal projections $\frac{1}{2}(p+s_{j}p)$
of $p$ to the hyperplane $H_{j}$. Moreover there is no vertex of $P_{h}$ contained in $H_{j}$, because $h^{c}$
is strictly convex on $\Delta^{c}$.

The function $h_{K}$ associated to a  polyhedron $K$  has always finite values if and only if $K$ is compact.
Moreover, there is a decomposition in convex cones of the convex set $\{n\in N_{\textbf{R}}:
h_{K}(n)\in\textbf{R}\}$ such that there is an one-to-one correspondence between the cones of such decomposition
and the faces of $K$. Hence there is an one-to-one correspondence between the 1-dimensional cones of such
decomposition and the semi-spaces that define $K$. Given a such cone $\tau$ the associated semi-space is $\{m\in
M_{\textbf{R}}: m(\varrho(\tau))\geq h_{K}(\varrho(\tau))\}$.

First of all we will show that $P_{h}\cap C^{+}=Q_{h}\cap C^{+}$.  It is sufficient to show that $Q_{h}\cap
C^{+}\cap M_{\textbf{Q}}\subseteq P_{h}$ because $P_{h}$ is closed. Given any $m\in Q_{h}\cap C^{+}\cap
M_{\textbf{Q}}$ and any $\tau\in\Delta^{c}(1)$ we have to show that $m(\varrho(\tau))\geq h^{c}(\varrho(\tau))$.
Because of the symmetry of $\Delta^{c}$, there are $w\in W^{1}$ and $\tau'\in\Delta(1)$ such that
$\varrho(\tau)=w\cdot \varrho(\tau')$. Observe that $w^{-1}\cdot m -m$ is a linear combination $\sum c_{i}f_{i}$
of the $f_{i}$ with positive coefficients, so $m(\varrho(\tau))=m(w\cdot \varrho(\tau'))= (w^{-1}\cdot
m)(\varrho(\tau'))= m(\varrho(\tau'))+\sum c_{i}f_{i}(\varrho(\tau')) \geq m(\varrho(\tau'))\geq
h(\varrho(\tau'))=h^{c}(\varrho(\tau))$. We now can conclude the proof. The decomposition in cones of
$N_{\textbf{R}}$ associated to $h_{P_{h}\cap C^{+}}$ has 1-dimensional cones
$\{\sigma(\check{g}_{1}),...,\sigma(\check{g}_{l})\} \cup \Delta(1)$. $h_{P_{h}\cap C^{+}}$ has  finite values
on all $N_{\textbf{R}}$, it is equal to $h$ on $|\Delta|$ and vanishes on the vectors $\check{g}_{1},...,
\check{g}_{l}$. The function associated to $\sigma(f_{1},...,f_{l})$ vanishes on $|\Delta|$ and has value
$-\infty$ on the complementary set. Thus their sum is the function $h$ associated to $Q_{h}$ and the proposition
follows by the fact that $h_{Q}+h_{Q'}=h_{Q+Q'}$ for each polyhedrons $Q$ and $Q'$. $\square$

We can prove a  stronger statement on the "rational" points of $Q_{h}$ and $P_{h}$.

\begin{prop}\label{red otor l1b}
Let $L_{h}$ be an ample spherical line bundle on $Y$, then $Q_{h}\cap (v_{h}+M)=P_{h}\cap C^{+}\cap
(v_{h}+M)+\sum_{i=1}^{l}\textbf{Z}^{+}f_{i}$.
\end{prop}

{\em Remark.} Observe that $H^{0}(Z,L|Z)$ is a
$O_{Z^{c}}(Z^{c})$-module through the restriction map
$O_{Z^{c}}(Z^{c})\rightarrow O_{Z}(Z)$ and $H^{0}(Z^{c},L|Z^{c})$ is
a $O_{Z^{c}}(Z^{c})$-submodule of $H^{0}(Z,L|Z)$. This proposition
implies that $H^{0}(Z,L|Z)$ is generated by  $H^{0}(Z^{c},L|Z^{c})$
as an $O_{Z}(Z)$-module.

{\em Proof.}  Observe that $f_{j}$ is orthogonal to $H_{j}$ and let
$\widetilde{f}_{i}=\frac{1}{2}(f_{i}+s_{i}f_{j})$ for each $i\neq j$. Thus  $\widetilde{f}_{i}\in H_{j}$ for
each $i\neq j$ and $\{\widetilde{f}_{i}\}_{i\neq j}$ is a basis of $H_{j}$. $-f_{i}$ and $-f_{j}$ are distinct
simple restricted roots, so $\widetilde{f}_{i}=f_{i}+d_{i}f_{j}$ for a suitable positive integer $d_{i}$. We
have the following easy consequence of the proposition~\ref{red otor l1a}.

\begin{lem}
$Q_{h}\cap H_{j}=P_{h}\cap H_{j}\cap C^{+}+\bigoplus_{i\neq j} \textbf{R}^{+}\tilde{f}_{i}$.
\end{lem}

{\em Proof.} Given $p=p'+\sum r_{i}f_{i}\in Q_{h}\cap H_{j}$ with $p'\in P_{h} \cap C^{+}$ and  $r_{i}$ positive
constants, we have $p=\frac{1}{2}(p'+s_{j}p')+\sum r_{i}\frac{1}{2}(f_{j}+s_{j}f_{i})$. $\square$

Let $R_{j}=\{p+af_{i}\ |\ p\in Q_{h}\cap H_{j}$ and $ -1/2\leq a\leq
1/2\}$. First of all we want to describe the conditions for a point
$m\in M_{\textbf{R}}$ to belong to $R_{j}$. Fixed any $j$, we define
another basis $u_{1},...,u_{l}$  of $M_{\textbf{R}}$ such that
$u_{j}=f_{j}$ and $u_{i}=g_{i}$ if $i\neq j$. The conditions for a
point $p=\sum y_{i}u_{i}$ to belong to $Q_{h}\cap H_{j}$ are
$y_{j}=0$ plus conditions of the form $\sum_{i\neq j}n_{i}y_{i}\geq
n$. Thus the conditions for a point $p=\sum y_{i}u_{i}$ to belong to
$R_{j}$ are the inequalities of the form $\sum_{i\neq
j}n_{i}y_{i}\geq n$ that define $Q_{h}\cap H_{j}$ plus the
inequalities $-1/2\leq y_{j} \leq1/2$. The following fundamental
lemma is the unique part of the proof in which we will use the
strictly convexity of $h^{c}$. (Remember that $h^{c}$ is strictly
convex on $\Delta^{c}$ if and only if $L_{h}$ is ample).

\begin{lem}\label{red otor l2}
$R_{j}$ is contained in $Q$ for each $j$.
\end{lem}

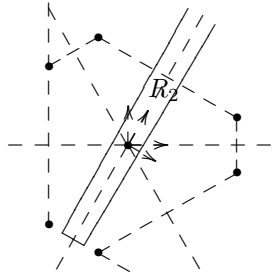
\[\begin{xy} <3pt, 0pt>:(0,0)*{\scriptstyle\bullet}="o";
(-10,-10)*{\scriptstyle\bullet}="p1"; (-10,10)*{\scriptstyle\bullet}="p2";
(-3.7648,13.6)*{\scriptstyle\bullet}="p3";"p3"+(17.4932,-10.1)*{\scriptstyle\bullet}="p4";
 (-3.7648,-13.6)*{\scriptstyle\bullet}="p8"; "p8"+(17.4932,10.1)*{\scriptstyle\bullet}="p7";
"o"+0;"o"+(5,0)**\dir{--}?>*\dir{>}; "o"+0;"o"+(18,0)**\dir{--}; "o"+0;"o"+(-15,0)**\dir{--};
"o"+0;"o"+(2.5,4.3333)**\dir{--}?>*\dir{>}; "o"+0;"o"+(9.477,16.416999)**\dir{--};
"o"+0;"o"+(-9,-15.588)**\dir{--};
 "o"+0;"o"+(-10,17.320)**\dir{--}; "o"+0;"o"+(9,-15.588)**\dir{--};
"o"+0;"o"+(3.4640,-2)**\dir{--}?>*\dir{>};"o"+0;"o"+(0,5)**\dir{--}?>*\dir{>};
 "p1"+0;"p2"+(0,8)**\dir{--}; "p2"+0;"p3"+0**\dir{--}; "p3"+0;"p4"+0**\dir{--}; "p4"+0;"p7"+0**\dir{--};
   "p7"+0;"p8"+0**\dir{--};"p8"+0,"p8"+(3.9905,-2.304)**\dir{--};
   "p8"+0;"p8"+(-1.7320,1)**\dir{--}; "p1"+(1.7320,-1);"p8"+(-1.7320,1)**\dir{-}; "p1"+0;"p1"+(1.7320,-1)**\dir{--};
  "p8"+(-1.7320,1);"p8"+(14.79708,28.8684)**\dir{-}; "p1"+(1.7320,-1);"p1"+(17.60292,27.2484)**\dir{-};
"o"+(4.5,7)*{R_{2}};
\end{xy}\]

{\em Proof.} Observe that it is sufficient to show that $P_{h}\cap H_{j}\cap C^{+}+[-1/2,1/2]f_{j}\subset Q_{h}$
because of the previous lemma. Because of the convexity of $Q_{h}$ it is sufficient to show that $Q_{h}$
contains the points $p'\pm (1/2)f_{j}$ for each vertex $p'$ of $P_{h}\cap H_{j}$.

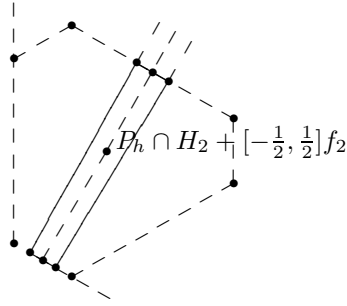
\[\begin{xy} <3.5pt, 0pt>:(0,0)*{\scriptstyle\bullet}="o";
(-10,-10)*{\scriptstyle\bullet}="p1"; (-10,10)*{\scriptstyle\bullet}="p2";
(-3.7648,13.6)*{\scriptstyle\bullet}="p3";"p3"+(17.4932,-10.1)*{\scriptstyle\bullet}="p4";
 (-3.7648,-13.6)*{\scriptstyle\bullet}="p8"; "p8"+(17.4932,10.1)*{\scriptstyle\bullet}="p7";
(-6.8824,-11.8)*{\scriptstyle\bullet}="q1";"p3"+(8.7466,-5.05)*{\scriptstyle\bullet}="q2";
"p8"+(-1.7320,1)*{\scriptstyle\bullet}="q1-";
"p1"+(1.7320,-1)*{\scriptstyle\bullet}="q1+";"q2"+(1.7320,-1)*{\scriptstyle\bullet}="q2-";
"q2"+(-1.7320,1)*{\scriptstyle\bullet}="q2+";
"o"+0;"o"+(7.5,12.999)**\dir{--};
"p8"+(-3.1176,1.8);"o"+0**\dir{--};"q1+"+0;"q1-"+0**\dir{-};"q2+"+0;"q2-"+0**\dir{-}; "q1+"+0;"q2+"+0**\dir{-};
"q2-"+0;"q1-"+0**\dir{-};"p8"+0,"p8"+(4.1568,-2.4)**\dir{--}; "p1"+0;"p2"+(0,6)**\dir{--};
"p2"+0;"p3"+0**\dir{--}; "p3"+0;"p4"+0**\dir{--}; "p4"+0;"p7"+0**\dir{--};
 "p7"+0;"p8"+0**\dir{--};"q1-"+0;"p8"+0**\dir{--};
"q1+"+0;"p1"+0**\dir{--};
    "q2-"+0;"q2-"+(2.5,4.3330)**\dir{--}; "q2+"+(0,0);"q2+"+(2.5,4.3330)**\dir{--};
 "o"+(13.5,1)*{P_{h}\cap H_{2}+[-\frac{1}{2},\frac{1}{2}]f_{2}};
\end{xy}\]

Observe that the vertices of $P_{h}\cap H_{j}$ are orthogonal projections to $H_{j}$ of suitable vertices of
$P_{h}$. Indeed let $p'$ a vertex of  $P_{h}\cap H_{j}$ and let $p$ be the endpoint different by $p'$ of the
segment intersection of  $P_{h}$ with the semi-line outgoing from $p'$ and parallel to $f_{j}$. If $p$ is not a
vertex of $P_{h}$ then $p$ is an interior point of a segment $I$ contained in $P_{h}$. Thus $p'$ is an interior
point of the projection of $I$ to $H_{j}$   and this segment is contained in $P_{h}$ by the symmetry of $P_{h}$,
a contradiction.

If $q'+af_{j}$ with $q'\in Q_{h}\cap H_{j}$ belongs to $M+v_{h}$, then $s_{j}(q'+af_{j})=q'-af_{j}$, so $2a\in
\textbf{Z}$. Moreover if $q$ is a vertex of $P_{h}\cap H_{j}$, then there is a constant $a$ such that $q+af_{j}$
is a vertex of $P_{h}$, so it is sufficient to show that the intersection of $P_{h}$ with the line parallel to
$f_{j}$ and passing through any vertex of $P_{h}\cap H_{j}$ is not a point.   If there is a vertex $p$ of
$P_{h}\cap H_{j}$ without such property, then $p$ is vertex of $P_{h}$  belonging to $H_{j}$, a contradiction.
$\square$

Now, we can conclude the proof of the proposition~\ref{red otor l1b} (look to the following figure). Let  $p$ be
a point contained in $Q_{h}\cap (M+v_{h})$ and suppose that $p=\sum x_{i}f_{i}= \sum \dot{x}_{i}g_{i}$. If
$\dot{x}_{i}\leq 0$ for each $i$, then $p\in P_{h}\cap C^{+}$. Otherwise there is an index $j$ such that
$\dot{x}_{j}>0$. We know that $p=p'+\sum a_{i}f_{i}$ where $p'\in P\cap C^{+}$ and the  $a_{i}$ are positive
constants. If $\dot{x}_{j}\geq2$ then $a_{j}\geq1$.  Thus the point
$p-[a_{j}]f_{j}=p'+(a_{j}-[a_{j}])f_{j}+\sum_{i\neq j} a_{i}f_{i}$ belongs to $Q_{h}\cap (M+v_{h})$ and it has
$j$-th coordinate with respect to $\{g_{1},...,g_{l}\}$ strictly less than 2 ($[a_{j}]$ is the integral part of
$a_{j}$). Moreover, this coordinate can be at most 1 because $p-[a_{j}]f_{j}$ is a weight. We can suppose that
it is exactly 1, so $p-[a_{j}]f_{j}-(1/2)f_{j}$ belongs to $Q_{h}\cap H_{j}$ and it is  the projection of
$p-[a_{j}]f_{j}$ to $H_{j}$. Thus $p-[a_{j}]f_{j}$ belong to $R_{j}$, so also $p-[a_{j}]f_{j}-f_{j}$ belongs to
$R_{j}$ and its $j$-th coordinate with respect to $\{g_{1},...,g_{l}\}$ is negative. Moreover
$p-(p-[a_{j}]f_{j}-f_{j})$ is a linear combination of the $f_{i}$ with positive integral coefficients. If there
is an index $k$ such that $p-[a_{j}]f_{j}-f_{j}$ has negative $k$-th coordinate with respect to
$\{g_{1},...,g_{l}\}$, then we reiterate the process. The process has to end in a finite number of steps because
$Q_{h}$ is contained in the semi-space $\{\sum x_{i}\geq h(e_{1}+...+e_{l})\}$. $\square$

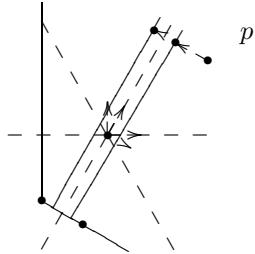
\[\begin{xy} <2.5pt, 0pt>:(0,0)*{\scriptstyle\bullet}="o";
(-10,-10)*{\scriptstyle\bullet}="p1"; (-10,10)*{\scriptstyle}="p2";
  (3.7648,-13.6)*{\scriptstyle}="p7";(-3.7648,-13.6)*{\scriptstyle\bullet}="p8";
 (15.1125,11.25)*{\scriptstyle\bullet}="x"; (10.2553,13.961)*{\scriptstyle\bullet}="y";
 (6.99383,15.886)*{\scriptstyle\bullet}="z";
"x"+0;"y"+0**\dir{--}?>*\dir{>}; "y"+0;"z"+0**\dir{--}?>*\dir{>};"o"+0;"o"+(5,0)**\dir{--}?>*\dir{>};
"o"+0;"o"+(15,0)**\dir{--}; "o"+0;"o"+(-15,0)**\dir{--}; "o"+0;"o"+(2.5,4.3330)**\dir{--}?>*\dir{>};
"o"+0;"o"+(9.5625,16.5645)**\dir{--}; "o"+0;"o"+(-10,-17.320)**\dir{--};
 "o"+0;"o"+(-10,17.320)**\dir{--}; "o"+0;"o"+(10,-17.320)**\dir{--};
            "o"+0;"o"+(3.4640,-2)**\dir{--}?>*\dir{>};"o"+0;"o"+(0,5)**\dir{--}?>*\dir{>};
 "p1"+0;"p2"+(0,10)**\dir{-};
   "p8"+0;"p1"+0**\dir{-};"p8"+0;"p8"+(6.928,-4)**\dir{-};
  "p8"+(-1.7320,1);"p8"+(15.0711,29.403)**\dir{-};
   "p1"+(1.7320,-1);"p1"+(17.9289,27.753)**\dir{-};
"o"+(21,15)*{p};
\end{xy}\]

Now we can prove the most important theorem of this work.

\begin{thm}\label{red otor}
Let  $L_{h}$ and $L_{k}$ be two spherical ample line bundles on $Y$. Then $m_{h,k}$ is surjective if and only if
$m^{c}_{h,k}$ is surjective.
\end{thm}

One can easily show that the theorem is equivalent to the following more combinatorial statement:

\[Q_{h}\cap (v_{h}+M) +Q_{k}\cap (v_{k}+M) =Q_{h+k}\cap (v_{h+k}+M) \]
if and only if
\[P_{h}\cap (v_{h}+M) +P_{k}\cap (v_{k}+M) =P_{h+k}\cap (v_{h+k}+M). \]

{\em Proof.} The sufficiency of the condition  is  easy. Given a point $p\in Q_{h+k}\cap (M+v_{h+k})$ we know
that $p=p'+\sum c_{i}f_{i}$ where $p'\in P_{h+k}\cap C^{+} \cap (M+v_{h+k})$ and the $c_{i}$ are positive
integers. Moreover  there are $p_{h}\in P_{h}\cap (M+v_{h})$ and $p_{k}\in P_{k}\cap (M+v_{k})$ such that
$p'=p_{h}+p_{k}$. Thus  $p=(p_{h}+\sum c_{i}f_{i})+p_{k}$  and $p_{h}+\sum c_{i}f_{i}$ belongs to  $Q_{h}\cap
(M+v_{h})$.

Suppose now that $Q_{h}\cap (v_{h}+M) +Q_{k}\cap (v_{k}+M) =Q_{h+k}\cap (v_{h+k}+M)$. Let  $m=\sum
z_{i}f_{i}=\sum \dot{z}_{i}g_{i}$ be a point in $P_{h+k}\cap (M+v_{h+k})$. We can suppose that $m$ belongs  to
$C^{+}$ by the symmetry of the polytopes $P_{h}$ and $P_{k}$. By hypothesis there are two points $p'_{0}\in
Q_{h}\cap (M+v_{h})$ and $q_{0}'\in Q_{k}\cap (M+v_{k})$ such that $p_{0}'+q_{0}'=m$. First, we will show that
we can choose $p'_{0}$ and $q'_{0}$ such that $p_{0}'$ belongs to $P_{h}$. Indeed we know that $p'_{0}=p_{0}+w$
where $p'_{0}\in P_{h}\cap C^{+}\cap (M+v_{h})$ and $w\in \bigoplus \textbf{Z}^{+}f_{i}$, so $m=p_{0}+q_{0}$
where $q_{0}:=q'_{0}+w$ belongs to $ Q_{k}\cap (M+v_{k})$.

Proceeding as in the proposition~\ref{red otor l1b}, we can define a sequence of pairs of points
$\{(p_{i},q_{i})\}_{i=0,...,r}$ with the following properties: 1) $p_{i}\in Q_{h}\cap (M+v_{h})$ for each $i$;
2)  $q_{i}\in Q_{k}\cap (M+v_{h})$ for each $i$; 3) $m=p_{i}+q_{i}$ for each $i$; 4) $(p_{0},q_{0})$ is as
before; 5) $(p_{i+1},q_{i+1})=(p_{i}+f_{j_{i}},q_{i}-f_{j_{i}})$ for a suitable  $j_{i}$ and 6) $q_{r}\in
P_{k}$. Indeed we can define the $\{q_{i}\}$ as in the proposition~\ref{red otor l1b} and then we set
$p_{i}=m-q_{i}$. Now it is sufficient to show by induction that we can choose the indices $j_{i}$ so that
$p_{i}$ belongs to $P_{h}$ for each $i$. We known that $p_{0}\in P_{h}$. Now suppose that $p_{n}$ belongs to
$P_{h}$ by inductive hypothesis. Suppose that $p_{n}=\sum x_{i}f_{i}=\sum \dot{x}_{i}g_{i}$ and $q_{n}=\sum
y_{i}f_{i}= \sum \dot{y}_{i}g_{i}$. If $q_{n}\in P_{k}$ we define  $r=n$ and there is nothing to prove.
Otherwise there is an index $j_{n}$ such that $\dot{y}_{j_{n}}>0$ and it is sufficient to prove that
$p_{n}+f_{j_{n}}$ belongs to $P_{h}$. Observe that $-\dot{x}_{j_{n}}> 0$, so $-\dot{x}_{j_{n}}\geq 1$ because it
is an integer. Moreover $s_{j_{n}}p_{n} =p_{n}-\dot{x}_{j_{n}} f_{j_{n}}$ belongs to $P_{h}$. Thus $P_{h}$
contains $p_{n}+f_{j_{n}}$ because it is convex. Thus we can choose $p_{n+1}=p_{n}+f_{j_{n}}$. $\square$

{\em Remark.} 1) The previous theorem is valid with the weaker hypotheses that $h$, $k$ are convex and that
$h|\sigma$, $k|\sigma$ are regular spherical weights for each $\sigma\in\Delta(l)$. Indeed these hypotheses
implies that no vertex of $P_{h}$ (respectively of $P_{k}$) is contained in a Weyl wall.

2) Suppose that $h=k$ is convex and that $h|\sigma$ is a  regular spherical weights for each
$\sigma\in\Delta(l)$. In this case one can show that $L_{h}|Z$ is the pullback of  an ample line bundle on a
possibly singular toric variety $Z'$ over $\textbf{A}^{l}$. This suggests to consider only ample line bundles.

\section{Line bundles on an exceptional complete symmetric variety}

Let $Y$ be an exceptional complete symmetric variety, let $Z$ be the associated open toric variety and let
$\Delta$ be the fan of $Z$. Given an ample  spherical line bundle $L_{h}$ over $Y$, we know that that the
multiplication $M_{h,h}$ of sections on $Y$ is surjective if and only if the multiplication $m_{h,h}$ of
sections on $Z$ is surjective. In this section we want to generalize this fact to the not spherical line
bundles. Remember that $Pic(X)$ is generated by the spherical weights and by the fundamental weights
$\omega_{\alpha_{1}},...,\omega_{\alpha_{s}}$ corresponding to the exceptional roots
$\alpha_{1},...,\alpha_{s}$.

\begin{prop} Let $L_{h'}$ be an ample line bundle on $Y$ such that $M_{h',h'}$ is surjective and let
$a_{1},...,a_{l}$ be positive integers. If we define $h=h'+\sum a_{i}\omega_{\alpha_{i}}$ then the product
$M_{h,h}$ of sections of $L_{h}$ over $Y$ is surjective.
\end{prop}

{\em Proof.} Observe that $L_{h}$ is an ample bundle on $Y$. We will prove the proposition by induction on $\sum
a_{i}$. $M_{h,h}$ is trivially surjective if $\sum a_{i}=0$. We need a lemma on the maps
$M_{h,\omega_{\alpha_{i}}}$.

\begin{lem}
Let $L_{h}$ be an ample line bundle on $Y$ and let $\omega\in\{\omega_{\alpha_{1}},...,\omega_{\alpha_{s}}\}$.
Then $M_{h,\omega}$ is surjective.
\end{lem}

{\em Proof.} In the following  $V_{\lambda}^{*}$ is the unique subrepresentation of $H^{0}(Y,L_{\lambda})$ which
contains a  lowest weight vector $v_{\lambda}$ of weight $-\lambda$. We have
$H^{0}(Y,L_{h})=\bigoplus_{\lambda\in \Pi(Y,h)}s^{h-\lambda}V^{*}_{\lambda}$, $H^{0}(Y,L_{h+\omega}) =
\bigoplus_{\lambda\in \Pi(Y,h)}s^{h-\lambda}V^{*}_{\omega+\lambda}$ and $H^{0}(Y,L_{\omega})=V^{*}_{\omega}$.
The lemma is implied by the fact that, for each $\lambda\in \Pi(Y,h)$,
$M_{h,\omega}(s^{h-\lambda}v_{\lambda}\otimes v_{\omega})$ is a lowest weight vector of weight
$-\lambda-\omega$. $\square$

We now go back to the proposition.  Let $j$ be an index such that
$a_{j}>0$ and define $\widetilde{h}=h-\omega_{\alpha_{j}}$. We have
the following commutative diagram

\[\xymatrix{ H^{0}(Y,L_{\widetilde{h}})\otimes H^{0}(Y,L_{\widetilde{h}})\otimes H^{0}(Y,L_{\omega_{\alpha_{j}}})
\otimes H^{0}(Y,L_{\omega_{\alpha_{j}}})\ar[r]\ar[d]^{m_{1}}
 & H^{0}(Y,L_{h})\otimes H^{0}(Y,L_{h})\ar[dd]^{M_{h,h}} \\
H^{0}(Y,L_{2\widetilde{h}})\otimes
H^{0}(Y,L_{\omega_{\alpha_{j}}})\otimes
H^{0}(Y,L_{\omega_{\alpha_{j}}})\ar[d]^{m_{2}}
&  \\
H^{0}(Y,L_{2\widetilde{h}+\omega_{\alpha_{j}}})\otimes
H^{0}(Y,L_{\omega_{\alpha_{j}}})\ar[r]^{M_{2\widetilde{h}+w_{\alpha_{j}},w_{\alpha_{j}}}}
& H^{0}(Y,L_{2h }).}\] $m_{1}$ is surjective by induction, $m_{2}$
and $M_{2\widetilde{h}+w_{\alpha_{j}},w_{\alpha_{j}}}$ are
surjective because of the previous lemma, so $M_{h,h}$ is
surjective. $\square$

\begin{thm}
Let $L_{h}$ be an ample line bundle on $Y$. If $m_{h,h}$ is surjective then $M_{h,h}$ is surjective.
\end{thm}

{\em Proof.} We know that, up to exchange $\alpha_{i}$ with $\overline{\theta}(\alpha_{i})$ for some $i$ in
$\{1,...,l\}$, there are positive integers $a_{1},...,a_{l}$ such that the line bundle  $L_{h'}$, with
$h'=h-\sum a_{i}w_{i}$, is spherical and ample. The restriction of $L_{h}$ to $Z$ is isomorphic to the
restriction of $L_{h'}$ to $Z$, so $m_{h',h'}$ is surjective. Thus $M_{h',h'}$ is surjective because of the
theorem~\ref{red otor}. Hence $M_{h,h}$ is surjective by the previous proposition. $\square$

\section{Open projectively normal toric varieties}

Now we want to describe some families of open toric varieties such that, if $L_{h}$  is an ample line bundle  on
a such variety, then the product $m_{h,h}$ of sections is surjective. One family is formed by all the varieties
of dimension 2 proper over $\textbf{A}^{2}$. Moreover we will find an infinite  number of varieties that have
such property for every given dimension. In some cases we will prove that, given any two ample line bundles
$L_{h}$ and $L_{k}$ on a fixed variety, then the product $m_{h,k}$ is surjective. {\em In the following we will
identify $M$ with $\textbf{Z}^{l}$}.

\subsection{Blow-ups of $\textbf{A}^{l}$}

Now we study the class of varieties that are  blow-ups of
$\textbf{A}^{l}$ along an irreducible stable closed subvariety.

\begin{prop}\label{ex 1}
Let  $Z$ be the blow-up of $\textbf{A}^{l}$ along the irreducible
stable closed subvariety associated to $\sigma(e_{1},...,e_{r})$.
Let  $L_{h}$ and $L_{k}$ be  two line bundles generated by global
sections on $Z$, then the product of sections $m_{h,k}$ is
surjective.
\end{prop}

The inequalities for $Q_{h}$ are $z_{i}\geq a_{i}$ for each  $i=1,..,l$ and $z_{1}+...+z_{r}\geq b$. The
inequalities for $Q_{k}$ are $z_{i}\geq c_{i}$ for each $i=1,..,l$ and $z_{1}+...+z_{r}\geq d$. Here $a_{i}$,
$b$, $c_{i}$ and $d$ are suitable integers. Let  $m$ be any point in $ Q_{h+k}\cap  M$, then there are
$\tilde{m}_{1}\in Q_{h}$ and $\tilde{m}_{2}\in Q_{k}$ such that $\tilde{m}_{1}+\tilde{m}_{2}=m$, but they may
have not integral coordinates. We want to translate $\tilde{m}_{1}$ with respect to a "little" vector $v$ so
that $\tilde{m}_{1}+v$ and $\tilde{m}_{2}-v$ will belong respectively to $Q_{h}\cap\textbf{Z}^{l}$ and to
$Q_{k}\cap\textbf{Z}^{l}$.  If $\tilde{m}_{1}=(x_{1},...,x_{l})$ then $x_{i}\geq a_{i}$. Let $[x_{i}]$ be the
integral part of $x_{i}$ and let $\epsilon_{i}=-[([x_{i}]-x_{i})]$ ($\epsilon_{i}$ is 0 if $x_{i}$ is an integer
and it is 1 otherwise). $[x_{i}]+\epsilon_{i} \geq[x_{i}]\geq a_{i}$ because the $a_{i}$ are integers. Likewise,
if $\tilde{m}_{2}=(y_{1},...,y_{l})$ then $[y_{i}]+\epsilon_{i} \geq[y_{i}] \geq d_{i}$ and
$\epsilon_{i}=-[([y_{i}]-y_{i})]$. If $([y_{1}],...,[y_{l}])$ belongs to $Q_{k}$, then we define
$m_{1}=([x_{1}]+\epsilon_{1},...,[x_{l}]+\epsilon_{l})$ and $m_{2}=([y_{1}],...,[y_{l}])$. Clearly these points
satisfy our requests. In the same way, if $([x_{1}],...,[x_{l}])$ belongs to $Q_{h}$, then we define
$m_{1}=([x_{1}],...,[x_{l}])$ and $m_{2}=([y_{1}]+\epsilon_{1},...,[y_{l}]+\epsilon_{l})$. Thus we can suppose
that $\sum_{i=1}^{r}[x_{i}]\leq b$ and $\sum_{i=1}^{r}[y_{i}]\leq d$.  We define $m_{1}=
([x_{1}]+\epsilon_{1},...,[x_{s}]+\epsilon_{s}, [x_{s+1}],...,[x_{l}])$ for an index $s$ lesser  than $ r$ and
such that $h(e_{1}+...+e_{r}) = \sum_{i=1}^{r}[x_{i}]+\sum_{i=1}^{s}\epsilon_{i}= m_{1}(e_{1}+...+e_{r})$. There
is a such  $s$ because $\sum_{i=1}^{r}[x_{i}]-b$ is a negative integer,
$\sum_{i=1}^{r}([x_{i}]+\epsilon_{i})\geq b+d -\sum_{i=1}^{r}([y_{i}]) \geq b$ and $\epsilon_{i}\in\{0,1\}$ for
each $i$. Moreover we define $m_{2}=m-m_{1}$.  To verify that $m_{2}\in Q_{k}$ it is sufficient to show that
$\sum_{i=1}^{s}[y_{i}]+\sum_{i=s+1}^{r}([y_{i}]+\epsilon_{i})\geq k(e_{1}+...+e_{r})$. This is implied by the
inequality $m_{2}(e_{1}+...+e_{r})=(m-m_{1})(e_{1}+...+e_{r})\geq (h+k)(e_{1}+...+e_{r})-h(e_{1}+...+e_{r})$.
$\square$

Now we study a similar family of varieties, but  we require that the two line bundles $L_{h}$ and $L_{k}$ are
equal.

\begin{cor}\label{ex 1b}
Let  $Z$ be the open toric variety obtained from $\textbf{A}^{l}$ through the sequence of blow-ups along the
subvarieties associated respectively to $\sigma(e_{1},e_{2})$, $\sigma(e_{2},e_{3})$,...,
$\sigma(e_{r-1},e_{r})$. Let  $L_{h}$ be any line bundle generated by global sections on $Z$, then the product
of sections $m_{h,h}$ is surjective.
\end{cor}

{\em Proof.} The inequalities for $Q_{h}$ are: $ z_{i}\geq a_{i}$ for each $i=1,..,l$ and $z_{i-1}+z_{i}\geq
b_{i}$ for each $i=2,...,r$, where the $a_{i}$ and the  $b_{i}$ are suitable integers. Let $m=(x_{1},...,x_{l})
\in Q_{2h}\cap M=2Q_{h}\cap M$. Observe that $m'=(x_{1}/2,...,x_{l}/2)\in Q_{h}$ and $m'+m'=m$. We define
$\epsilon_{i}=-[([x_{i}/2]-x_{i}/2)]$, $m_{1}=m'+ (\epsilon_{1},0,\epsilon_{3}, ...,\epsilon_{s}, 0,...,0)$ and
$m_{2}=m-m_{1}$ for a suitable  $s$. If $r$ is odd then we choose $s=r$, otherwise we define $s=r-1$. These
points belong to $Q_{h}\cap M$ because $[x_{i-1}/2]+[x_{i}/2]+(\epsilon_{i-1}+\epsilon_{i})/2\geq b_{i}$ for
each $i$. $\square$

\subsection{Open toric varieties of dimension 2 and a singular family in dimension 3}

Now we consider the family  of  smooth toric varieties  proper over $\textbf{A}^{2}$.

\begin{thm}\label{ex 2}
Let  $Z$ be any smooth toric variety  proper over $\textbf{A}^{2}$. Let $L_{h_{1}}$ and $L_{h_{2}}$ be two
linearized line bundles generated by global sections and suppose that $h_{1}$ and $h_{2}$ are strictly convex on
the same fan, then the product of sections $m_{h_{1},h_{2}}$ is surjective.
\end{thm}

The hypotheses mean that there is a variety $Z'$ and two ample line bundle $L'_{h}$ and $L'_{k}$ over $Z'$ such
that $L_{h}$ is the pullback of $L'_{h}$ and $L_{k}$ is the pullback of $L'_{k}$. We want to remark that $Z'$
may be singular.

{\em Proof. } Define a scalar product $(\ ,\, )$ such that $\{f_{1},f_{2}\}$ is a orthonormal basis. In this
proof, when we will say that a side $L$ of a polytope $P$ is orthogonal to a vector $v$, we will always suppose
that $(p,v)\geq0$ for each $p\in P$ (and $(x,v)=0 $ for each $x\in L$). Notice that a plane $H$ is the locus of
zeroes of $x_{1}e_{1}+x_{2}e_{2}\in N$ if and only if it is orthogonal to $x_{1}f_{1}+x_{2}f_{2}$.

Let $h_{3}=h_{1}+h_{2}$ and let  $\Delta$ be the fan of $Z$. It is obviously   sufficient to prove that
$Q_{h_{3}}\cap M=Q_{h_{1}}\cap M+Q_{h_{2}}\cap M$. We want to decompose each $Q_{h_{i}}$ in more simple
polyhedrons. More precisely we will decompose each $Q_{h_{i}}$ in two types of polyhedrons with vertices in $M$:
1) cones of form $p+\sigma(f_{1},f_{2})$ for a suitable point $p$ and 2) triangles. These triangles will have a
particular form, indeed we require that the fan associated to any such triangle  $\widetilde{T}$ has
1-dimensional cones generated respectively by $-e_{1},-e_{2}$ and by an element of $\sigma(e_{1},e_{2})$.

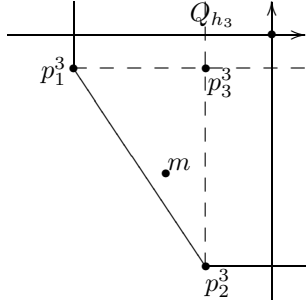
\[\begin{xy} <5pt, 0pt>:(30,47.5)*{\scriptstyle\bullet}="o";
(15,45)*{\scriptstyle\bullet}="p1";(25,30)*{\scriptstyle\bullet}="p2"; (22,37)*{\scriptstyle\bullet}="q";
(25,45)*{\scriptstyle\bullet}="p3"; "o"+0;"o"+(-20,0)**\dir{-}; "o"+0;"o"+(0,-20)**\dir{-};
"o"+0;"o"+(2.5,0)**\dir{-};"o"+0;"o"+(2.5,0)?>*\dir{>}; "o"+0;"o"+(0,2.5)**\dir{-};"o"+0;"o"+(0,2.5)?>*\dir{>};
"p1"+0;"p1"+(17.5,0)**\dir{--};"p1"+0;"p1"+(0,5)**\dir{-};
"p2"+0;"p2"+(7.5,0)**\dir{-};"p2"+0;"p2"+(0,20)**\dir{--}; "p1"+0;"p2"+0**\dir{-};
"p1"+(-1.5,-0.2)*{p_{1}^{3}};"p2"+(0.9,-1.3)*{p_{2}^{3}};"p1"+(11.1,-1.1)*{p_{3}^{3}};
"q"+(1,0.7)*{m};
"p1"+(10.5,4)*{Q_{h_{3}}};
\end{xy}\]

Let    $m=(x_{1},x_{2})$ be any point in  $Q_{h_{3}}\cap M$. If there is a vertex $p^{3}$ of $Q_{h_{3}}$ such
that $m$ is contained in the polyhedron $p^{3}+\sigma(f_{1},f_{2})$, then there is a maximal cone
$\sigma\in\Delta$ such that $p^{3}=h_{3}|\sigma=h_{1}|\sigma+h_{2}|\sigma$ and $m$ belongs to
$(h_{1}|\sigma+\sigma(f_{1},f_{2}))+ (h_{2}|\sigma+\sigma(f_{1},f_{2}))=p^{3}+\sigma(f_{1},f_{2})$ where
$h_{j}|\sigma+ \sigma(f_{1},f_{2})$ is the polyhedron associated to the pullback of a line bundle on
$\textbf{A}^{2}$. Otherwise there are two maximal cones $\sigma_{1}$ and $ \sigma_{2}$ with the following
properties. For each $j$ let $p_{1}^{j}=h_{j}|\sigma_{1}=(z_{1}^{1,j},z_{2}^{1,j})$, $p_{2}^{j}=
h_{j}|\sigma_{2} =(z_{1}^{2,j},z_{2}^{2,j})$, $p_{3}^{j} =(z_{1}^{2,j},z_{2}^{1,j})$ and define $T^{j}$ as the
triangle with vertices $p_{1}^{j}$, $p_{2}^{j}$ and $p_{3}^{j}$. Then $T^{3}=T^{1}+T^{2}$ and $T^{3}$ contains
$m$. Indeed we can choose $\sigma_{1}$ and $ \sigma_{2}$ such that: 1) $p_{1}^{3}$, $p^{3}_{2}$ are the two
vertices of a side of $Q_{h_{3}}$; 2) $z_{1}^{1,3}\leq x_{1}\leq z_{1}^{2,3}$. Observe \vspace{1 mm} that   $
x_{2}\leq z^{1,3}_{2}$ because otherwise 
$m$ belongs to $(z_{1}^{1,3},z_{2}^{1,3})+\sigma(f_{1},f_{2})$. Moreover the fans associated to these triangles
are equal to the fan with 1-dimensional \vspace{1 mm} cones $\sigma(-e_{1}), \sigma(-e_{2})$ and
$\sigma(a_{1}e_{1}+a_{2}e_{2})$ for suitable integers $a_{1}$ and $a_{2}$. Observe that $p_{1}^{j}$
(respectively $p_{2}^{j}$) is the  vertex of $T^{j}$ not contained in the side orthogonal to $-f_{1}$
(respectively to $-f_{2}$). Moreover we can suppose that $a_{1}f_{1}+a_{2}f_{2}\in\sigma(f_{1}+f_{2},f_{2})$ up
to exchange $f_{1}$ and $f_{2}$. In the pictures we consider the case in which $h_{1}=h_{2}$.

We want to prove that $(T^{1}+\sigma(f_{1},f_{2}))\cap M+
(T^{2}+\sigma(f_{1},f_{2}))\cap
M=(T^{1}+T^{2}+\sigma(f_{1},f_{2}))\cap M$ by induction on
$a_{1}+a_{2}$. If $a_{1}+a_{2}=2$ the previous equality is trivially
implied by the proposition~\ref{ex 1}. Otherwise we decompose
$T^{j}$ in two triangles by intersecting $T^{j}$ with the line
$r^{j}$ orthogonal to $f_{1}+f_{2}$ and passing for the vertex
$p_{1}^{j}$. Observe that $r^{j}$ intersects the side orthogonal to
$-f_{1}$ in the integral point
$q^{j}=(z^{2,3}_{1},z^{1,3}_{2}+z^{1,3}_{1}-z^{2,3}_{1})$. We have
decomposed $T^{j}$ in two triangles: 1) the triangle $T^{j}_{1}$
with sides orthogonal respectively to $-f_{1}$, $-f_{2}$,
$f_{1}+f_{2}$ and with vertices $p_{1}^{j}$, $p_{3}^{j}$, $q^{j}$;
2) the triangle $T^{j}_{0}$ with sides orthogonal respectively to
$-f_{1}$, $-f_{1}-f_{2}$, $a_{1}f_{1}+a_{2}f_{2}$ and with vertices
$p_{1}^{j}$, $q^{j}$, $p_{2}^{j}$.   Observe that $r^{1}_{i}$,
$r^{2}_{i}$ and $r^{3}_{i}$ are parallel, so the fans associated to
$T^{1}_{1}$, $T^{2}_{1}$ and $T^{3}_{1}$ (respectively to
$T^{1}_{0}$, $T^{2}_{0}$ and $T^{3}_{0}$) are equal.
$T^{3}_{1}=T^{1}_{1}+T^{2}_{1}$ because each $T^{j}_{1}$ share two
vertices  with $T^{j}$. In the same way we obtain that
$T^{3}_{0}=T^{1}_{0}+ T^{2}_{0}$. Notice that
$\sigma(f_{1}+f_{2},f_{2})\subset\sigma(f_{1},f_{2})$ so
$T_{0}^{j}+\sigma(f_{1}+f_{2},f_{2})\subset
T^{j}+\sigma(f_{1},f_{2})$ for each $j$. Moreover
$a_{1}f_{1}+a_{2}f_{2}= (a_{1}-a_{2})f_{1}+a_{2}(f_{1}+f_{2})$ and
$a_{1}<a_{1}+a_{2}$, so we can apply the inductive hypothesis.
$\square$

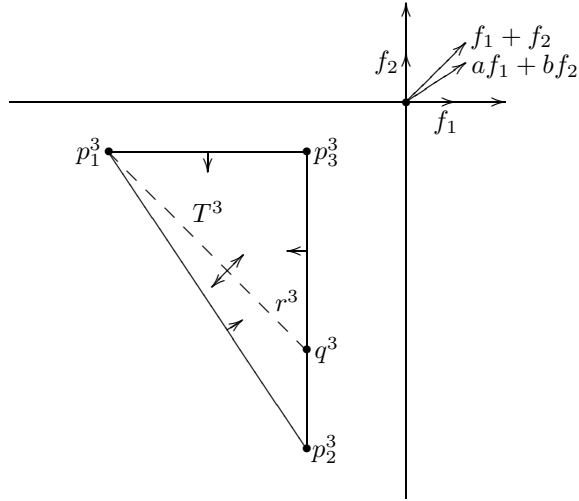
\[\begin{xy} <7.5pt, 0pt>:(30,47.5)*{\scriptstyle\bullet}="o";
(15,45)*{\scriptstyle\bullet}="p1";(25,30)*{\scriptstyle\bullet}="p2"; (25,35)*{\scriptstyle\bullet}="q";
(25,45)*{\scriptstyle\bullet}="p3"; "o"+0;"o"+(-20,0)**\dir{-}; "o"+0;"o"+(0,-20)**\dir{-};
"o"+0;"o"+(5,0)**\dir{-}?>*\dir{>};;"o"+0;"o"+(2.5,0)**\dir{-}?>*\dir{>};
"o"+0;"o"+(0,5)**\dir{-}?>*\dir{>};;"o"+0;"o"+(0,2.5)**\dir{-}?>*\dir{>}; "p1"+0;"p3"+0**\dir{-};
"p2"+0;"p3"+0**\dir{-};
 "p1"+0;"q"+0**\dir{--};"p1"+0;"p2"+0**\dir{-};
"p2"+(0,10);"p2"+(-1,10)**\dir{-};"p2"+(0,10);"p2"+(-1,10)?>*\dir{>};
"p1"+(5,0);"p1"+(5,-1)**\dir{-};"p1"+(5,0);"p1"+(5,-1)?>*\dir{>};
"p1"+(6,-6);"p1"+(6.8,-5.2)**\dir{-};"p1"+(6,-6);"p1"+(6.5,-5.5)?>*\dir{>};
"p1"+(6,-6);"p1"+(5.2,-6.8)**\dir{-};"p1"+(6,-6);"p1"+(5.5,-6.5)?>*\dir{>};
"p1"+(6,-9);"p1"+(6.75,-8.5)**\dir{-};"p1"+(6,-9);"p1"+(6.75,-8.5)?>*\dir{>};
"o"+0;"o"+(3,3)**\dir{-};"o"+0;"o"+(1.5,1.5)?>*\dir{>}; "o"+0;"o"+(3,2)**\dir{-};"o"+0;"o"+(1.5,1)?>*\dir{>};
"p1"+(-1,0)*{p_{1}^{3}};"p2"+(0.9,0)*{p_{2}^{3}};"p1"+(11,0)*{p_{3}^{3}}; "p1"+(5,-3)*{T^{3}};
"p2"+(-1,7.5)*{r^{3}};
"o"+(2,-1)*{f_{1}};"o"+(-1,2)*{f_{2}};"o"+(5.3,3.3)*{f_{1}+f_{2}};"o"+(6,1.8)*{af_{1}+bf_{2}};
"p1"+(11,-10)*{q^{3}};
\end{xy}\]

Now we consider a class of line bundles on varieties of dimension 3. This line bundles are the pullbacks of
ample lines bundles on varieties which are usually singular.

\begin{prop}\label{ex 2b}
Let  $\Delta$ be the fan with maximal cones $\sigma(e_{1},e_{2}, ae_{1}+ae_{2}+e_{3})$, $\sigma(e_{1},e_{3},
ae_{1}+ae_{2}+e_{3})$ and $\sigma(e_{2},e_{3}, ae_{1}+ae_{2}+e_{3})$ where $a$ is a strictly positive integer.
If $h$ is a strictly convex $\Delta$-linear function, then $Q_{2h}\cap M=Q_{h}\cap M+ Q_{h}\cap M$.
\end{prop}

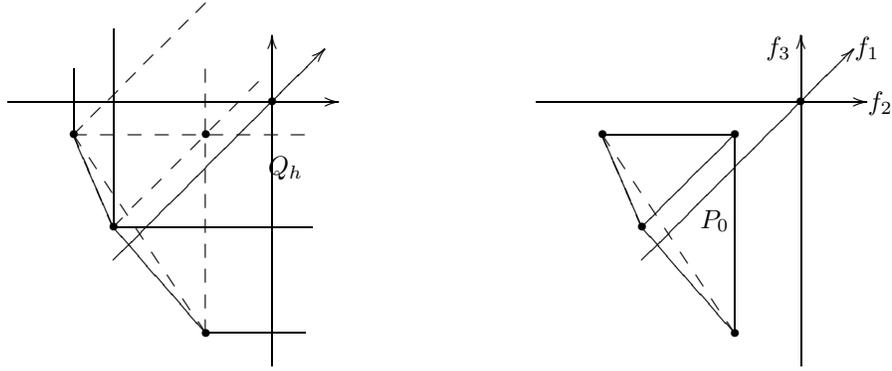
\[\begin{xy} <5pt, 0pt>:(30,47.5)*{\scriptstyle\bullet}="o";
(15,45)*{\scriptstyle\bullet}="p1";(25,30)*{\scriptstyle\bullet}="p2"; (18,38)*{\scriptstyle\bullet}="q";
(25,45)*{\scriptstyle\bullet}="p3"; "o"+0;"o"+(-20,0)**\dir{-}; "o"+0;"o"+(0,-20)**\dir{-};
"o"+0;"o"+(5,0)**\dir{-};"o"+0;"o"+(5,0)?>*\dir{>}; "o"+0;"o"+(0,5)**\dir{-}?>*\dir{>};
"o"+0;"o"+(4,4)**\dir{-}?>*\dir{>};"o"+0;"o"+(-12,-12)**\dir{-};
 "p1"+0;"p1"+(17.5,0)**\dir{--};"p1"+0;"p1"+(0,5)**\dir{-};
"p2"+0;"p2"+(7.5,0)**\dir{-};"p2"+0;"p2"+(0,20)**\dir{--};
"p1"+0;"q"+0**\dir{--};"p1"+0;"p2"+0**\dir{--};
"p1"+0;"q"+0**\dir{-};"p3"+(4,4);"q"+0**\dir{--};"p2"+0;"q"+0**\dir{-};"q"+0;"q"+(15,0)**\dir{-};
"q"+0;"q"+(0,15)**\dir{-};"p1"+0;"p1"+(10,10)**\dir{--};"o"+(1,-5)*{Q_{h}};
(70,47.5)*{\scriptstyle\bullet}="o'";
(55,45)*{\scriptstyle\bullet}="p'1";(65,30)*{\scriptstyle\bullet}="p'2";
(58,38)*{\scriptstyle\bullet}="q'";
(65,45)*{\scriptstyle\bullet}="p'3"; "o'"+0;"o'"+(-20,0)**\dir{-};
"o'"+0;"o'"+(0,-20)**\dir{-};
"o'"+0;"o'"+(5,0)**\dir{-};"o'"+0;"o'"+(5,0)?>*\dir{>};
"o'"+0;"o'"+(0,5)**\dir{-}?>*\dir{>};
"o'"+0;"o'"+(4,4)**\dir{-}?>*\dir{>};"o'"+0;"o'"+(-12,-12)**\dir{-};
"p'1"+0;"q'"+0**\dir{--};"p'1"+0;"p'2"+0**\dir{--};"p'1"+0;"p'3"+0**\dir{-};"p'2"+0;"p'3"+0**\dir{-};
"p'1"+0;"q'"+0**\dir{-};"p'3"+0;"q'"+0**\dir{-};"p'2"+0;"q'"+0**\dir{-};
"o'"+(-6.5,-9)*{P_{0}};
"o'"+(5,4)*{f_{1}};"o'"+(6,-0.1)*{f_{2}};"o'"+(-1.7,4)*{f_{3}};
\end{xy}\]

Remember that $h$ defines a line bundle generated by global sections on every toric variety proper over the
toric variety associated to $\Delta$. Moreover the toric variety associated to $\Delta$ is proper over
$\textbf{A}^{3}$ and it is smooth if and only if $a=1$. (Look to the figure for an example).

{\em Proof.} We want to proceed as in the previous theorem. We again define a scalar product such that
$\{f_{1},f_{2},f_{3}\}$ is an orthonormal basis and we can again decompose  $Q_{2h}$ in a simplex  $P_{0}$ and
some cones $p+\sigma(f_{1},f_{2},f_{3})$. Hence we can reduce ourselves to prove by induction on $a$ that
$(2P+\sigma(v_{1},v_{2},v_{3}))\cap M=(P+\sigma(v_{1},v_{2},v_{3}))\cap M +(P+\sigma(v_{1},v_{2},v_{3}))\cap M$
for each simplex $P$ with faces orthogonal respectively to $-v_{1},-v_{2},-v_{3}$ and to $av_{1}+av_{2}+v_{3}$,
where $\{v_{1},v_{2},v_{3}\}$ is an opportune base of $\textbf{Z}^{3}$ and $a$ is a positive integer (observe
that the polytope $P_{0}$ satisfies such requests). We can suppose, up to a translation, that the origin 0 is
the vertex of $P$ which does not belong to the face orthogonal to $av_{1}+av_{2}+v_{3}$. Let $(-b,0,0)$,
$(0,-b,0)$ and $(0,0,-c)$ be the other vertices of $P$. We have $c=ba$, so  $c\geq b$. If $a=1$ the equality is
verified because of the proposition~\ref{ex 1}. Otherwise we decompose $P$ intersecting it with the plane
orthogonal to $v_{1}+v_{2}+v_{3}$ and passing through the vertices $(-b,0,0)$ and $(0,-b,0)$. This plane
intersects the side of $P$ parallel to $\textbf{R}v_{3}$ in $(0,0,-b)$. We obtain two simplices with integral
vertices. The first one has faces orthogonal respectively to $-v_{1}$, $-v_{2}$, $-v_{3}$ and to
$v_{1}+v_{2}+v_{3}$. This simplex has vertices $(0,0,0)$, $(-b,0,0)$, $(0,-b,0)$ and $(0,0,-b)$. The second
simplex  has faces orthogonal respectively to $-v_{1}$, $-v_{2}$, $-v_{1}-v_{2}-v_{3}$ and to
$av_{1}+av_{2}+v_{3}$. This simplex $T$  has vertices $(-b,0,0)$,  $(0,-b,0)$, $(0,0,-b)$ and $(0,0,-c)$.
Observe that $\sigma(v_{1},v_{2}, v_{1}+v_{2}+v_{3})$ is contained in $\sigma(v_{1},v_{2},v_{3})$, so
$T+\sigma(v_{1},v_{2}, v_{1}+v_{2}+v_{3})$ is contained in $P+\sigma(v_{1},v_{2},v_{3})$. Moreover $T$ is a
simplex of the same type of $P$ and $av_{1}+av_{2}+v_{3}= (a-1)v_{1}+(a-1)v_{2}+(v_{1}+v_{2}+v_{3})$, i.e. the
coordinate with respect to the new basis are decreased. $\square$

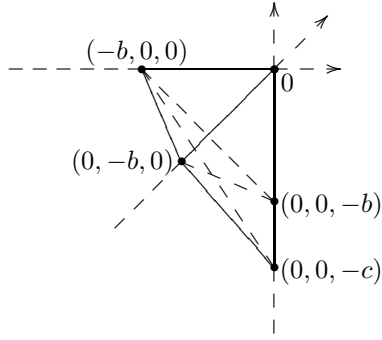
\[\begin{xy} <5pt, 0pt>:(30,47.5)*{\scriptstyle\bullet}="o";
(20,47.5)*{\scriptstyle\bullet}="p1";(30,32.5)*{\scriptstyle\bullet}="p2"; (23,40.5)*{\scriptstyle\bullet}="q";
(30,37.5)*{\scriptstyle\bullet}="t"; (30,47.5)*{\scriptstyle\bullet}="p3"; "p1"+0;"o"+(-20,0)**\dir{--};
"p2"+0;"o"+(0,-20)**\dir{--}; "o"+0;"o"+(5,0)**\dir{--};"o"+0;"o"+(5,0)?>*\dir{>};
"o"+0;"o"+(0,5)**\dir{--}?>*\dir{>}; "o"+0;"o"+(4,4)**\dir{--}?>*\dir{>};"q"+0;"o"+(-12,-12)**\dir{--};
"p1"+0;"p2"+0**\dir{--};"p1"+0;"p3"+0**\dir{-};"p2"+0;"p3"+0**\dir{-};
"p1"+0;"q"+0**\dir{-};"p3"+0;"q"+0**\dir{-};"p2"+0;"q"+0**\dir{-}; "t"+0;"q"+0**\dir{--};
"p1"+0;"t"+0**\dir{--};"p3"+(1,-1)*{0};
"p1"+(-0.4,1.3)*{(-b,0,0)};"p2"+(4.3,-0.3)*{(0,0,-c)};"q"+(-4.5,0.0)*{(0,-b,0)}; "t"+(4.3,-0.3)*{(0,0,-b)};
\end{xy}\]

\subsection{Stable subvarieties}

In some case we can reduce the study of the product of sections of
two ample line bundles to the study of the product of sections of
the restrictions of these line bundles to irreducible stable closed
subvarieties. First of all, Brion has proved the following
proposition.

\begin{prop}[See theorem 4 in \cite{BI}]\label{red svar 1}  Let $L$ be any  line bundle generated by global sections
on a smooth quasi-projective toric variety $Z$. Suppose that either
$Z$ is complete or it is proper over $\textbf{A}^{l}$. Then, given
two cones $\gamma\subset\gamma'$ in $\Delta$, the restriction
\[H^{0}(Z_{\gamma},L|Z_{\gamma})\longrightarrow H^{0}(Z_{\gamma'},L|Z_{\gamma'})\]
is surjective.
\end{prop}

\begin{prop}\label{red svar 3a} Let  $L_{h}$ and $L_{k}$ be two ample linearizated line bundles on
a smooth toric variety $Z$. Suppose that either $Z$ is complete or
it is proper over $\textbf{A}^{l}$. Let $\tau$ be a cone in
$\Delta(1)$ and let $s$ be a global section of $L_{h+k}$ which does
not vanish on $Z_{\tau}$. If $s|Z_{\tau}$ belongs to the image of
the product $m^{\tau}_{h,k}$ of  sections of the restrictions of
$L_{h}$ and $L_{k}$ to $Z_{\tau}$, then  $s$ belongs to the image of
the product $m_{h,k}$ of sections of $L_{h}$ and $L_{k}$.
\end{prop}

{\em Proof.} Since there is a basis of semi-invariant sections, we
can suppose that $s$ is a semi-invariant section  of weight  $\mu$,
so $\mu(\varrho(\tau))= (h+k)(\varrho(\tau))$.  Because of the
previous proposition, there are sections
$s'_{i}\in\Gamma(Z^{c},L_{h})$ and $s_{i}''\in\Gamma(Z^{c},L_{k})$
such that $m^{\tau}_{h,k}(\sum s'_{i}|Z_{\tau}^{c}\otimes
s_{i}''|Z_{\tau}^{c})$ is $s|Z^{c}_{\tau}$, so $m_{h,k}(\sum
s'_{i}\otimes s_{i}'')=s$. $\square$

\subsection{Two families of open toric varieties of dimension at least 3}

Now we want to show that there is an infinite number of open toric varieties of any fixed dimension (greater
than 2) such that the product of sections of any two ample line bundles is surjective. The principal instrument
in what follows  is the proposition~\ref{red svar 3a}. We will consider a very special class of varieties.
Indeed, given any ample line bundle $L$   on a variety of this family, then $H^{0}(Z,L)$ is generated as a
$O_{Z}(Z)$ module by the sections that do  not vanish on a suitable divisor.

\begin{prop}\label{ex 3-1}
Let $Z^{n}$ be the open toric variety obtained from $A^{l}$ through
the sequence of blow-ups along the irreducible stable subvarieties
associated respectively to $\sigma(e_{1},...,e_{l})$,
$\sigma(e_{1},...,e_{l-1},$ $ (\sum_{i=1}^{l-1}e_{i})+e_{l})$,
$\sigma(e_{1},...,e_{l-1},
2(\sum_{i=1}^{l-1}e_{i})+e_{l})$,...,$\sigma(e_{1},$ $...,$
$e_{l-1},$ $ i(\sum_{i=1}^{l-1}e_{i})
+e_{l})$,...,$\sigma(e_{1},...,e_{l-1},
(n-1)(\sum_{i=1}^{l-1}e_{i})+e_{l})$. Let $L_{h}$ and $L_{k}$ be any
two  line bundles generated by global sections on $Z^{n}$, then the
product of sections $m_{h,k}$ is surjective.
\end{prop}

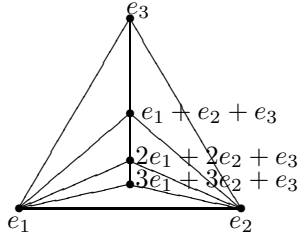
\[\begin{xy} <3pt, 0pt>: (0,0)*{\scriptstyle\bullet}="e1";
(28,0)*{\scriptstyle\bullet}="e2";(14,24)*{\scriptstyle\bullet}="e3"; (14,12)*{\scriptstyle\bullet}="v1";
(14,6)*{\scriptstyle\bullet}="v2"; (14,3)*{\scriptstyle\bullet}="v3";
"e1"+0;"e2"+0**\dir{-};"e1"+0;"e3"+0**\dir{-};"e2"+0;"e3"+0**\dir{-};
"v1"+0;"e1"+0**\dir{-};"v1"+0;"e2"+0**\dir{-};"v1"+0;"e3"+0**\dir{-};
"v2"+0;"e1"+0**\dir{-};"v2"+0;"e2"+0**\dir{-}; "v3"+0;"e1"+0**\dir{-};"v3"+0;"e2"+0**\dir{-};
"v2"+0;"v1"+0**\dir{-};"v2"+0;"v3"+0**\dir{-}; "e1"+(0,-2)*{e_{1}};"e2"+(0,-2)*{e_{2}};"e3"+(1,1)*{e_{3}};
"v1"+(10,0)*{e_{1}+e_{2}+e_{3}};"v2"+(11,0.5)*{2e_{1}+2e_{2}+e_{3}}; "v3"+(11,0.5)*{3e_{1}+3e_{2}+e_{3}};
\end{xy}\]

{\em Proof.} In the figure we have drawn the variety $Z^{3}$ with dimension 3. Observe that we have already
considered the case $n=1$ in   proposition~\ref{ex 1}, so we can suppose $n\geq 2$. Up to changing the
linearizations of the line bundles we can suppose that $h(e_{j})=k(e_{j})=0$ for each $j$. Observe that, if
$(Q_{h}\cap M)+(Q_{k}\cap M)$ contains a weight $p$, then it contains any weight $p+\sum a_{i}f_{i}$ where the
$a_{i}$ are positive integers. So we can  consider only the "minimal" weights. We now prove the propriety of $Z$
stated before the proposition.

\begin{cl}
Let $p$ be any weight in $Q_{h+k}\cap M$ and suppose that there is not a weight $p'$ in $Q_{h+k}\cap M$ such
that $p\in p'+\sigma(f_{1},...,f_{l})$. Then there is a cone $\tau\in\Delta(1)$ such that
$p(\varrho(\tau))=(h+k)(\varrho(\tau))$.
\end{cl}

{\em Proof.} Observe that $p(\varrho(\tau))=(h+k)(\varrho(\tau))$  means that any semi-invariant section of
weight $p$ does not vanish on the divisor $Z_{\tau}$. The hypotheses imply that $p-f_{l}$ does not belong to
$Q_{h+k}$. Hence there is an $i$ such that $(p-f_{l})(i(\sum_{i=1}^{l-1}e_{i})+e_{l})<
(h+k)(i(\sum_{i=1}^{l-1}e_{i})+e_{l})$, so $p(i(\sum_{i=1}^{l-1}e_{i})+e_{l})=
(h+k)(i(\sum_{i=1}^{l-1}e_{i})+e_{l})$. $\square$

Now it is sufficient to prove the surjectivity of the  product of sections of the restrictions of $L_{h}$ and
$L_{k}$ to the divisor $Z^{n}_{i}$ associated to $\sigma(i(\sum_{i=j}^{l-1}e_{j})+e_{l})$ for each $i=0,...,n$.
$Z_{0}$ is the blow-up of $A^{l-1}$ in the stable point;  $Z^{n}_{n}$ is the projective space of dimension $l-1$
while the other  $Z^{n}_{i}$ are isomorphic to the blow-up of  projective space in a point. Since $Z^{n}_{1}$
dominates $Z^{n}_{n}$ it is sufficient to study the product of sections of any two line bundles generated by
global sections on $Z^{n}_{1}$.

\begin{lem}\label{ex 3-1l}
Let $L_{h'}$ and $L_{k'}$ be any two line bundles on $Z^{n}_{1}$ generated by global sections. Then the
multiplication of sections is surjective.
\end{lem}

{\em Proof.} We can suppose that
$h'(\widetilde{e}_{i})=k'(\widetilde{e}_{i})=0$ for each $i$. In the
following we identify $\textbf{Z}^{l-1}$ with the character group
of the torus contained in $Z^{n}_{1}$. We proceed as in the proof of
the proposition~\ref{ex 1}. Given   any point $m$ in $Q_{h'+k'}$
with integral coordinates, there are $\tilde{m}_{1}\in Q_{h'}$ and
$\tilde{m}_{2}\in Q_{k'}$ such that $\tilde{m}_{1}+\tilde{m}_{2}=m$.
Now we want to simplify the notation. In particular we will be
evident that the problem does not depend on the dimension. Suppose
that $\tilde{m}_{1}=(x_{1},...,x_{l-1})$,
$\tilde{m}_{2}=(y_{1},...,y_{l-1})$ and $m=(z_{1},...,z_{l-1})$. Let
$[x_{i}]$ be the integral part of $x_{i}$ and let
$\epsilon_{i}=-[([x_{i}]-x_{i})]$. We define
$m_{1}=([x_{1}]+\epsilon_{1},...,[x_{\bar{r}}]+\epsilon_{\bar{r}},[x_{\bar{r}+1}],...,[x_{l-1}])$
and $m_{2} =([y_{1}],...,$ $[y_{\bar{r}}],$ $
[y_{\bar{r}+1}]+\epsilon_{\bar{r}+1},...,[y_{l-1}]+\epsilon_{l-1})$
for a suitable $\bar{r}$.  Let $t=\sum_{i=1}^{l-1} \epsilon_{i}$,
$r=\sum_{i=1}^{\bar{r}} \epsilon_{i}$,
$[x]=\sum_{i=1}^{l-1}[x_{i}]$, $x=\sum_{i=1}^{l-1}x_{i}$,
$[y]=\sum_{i=1}^{l-1}[y_{i}]$ and $y=\sum_{i=1}^{l-1}y_{i}$. Let
$a=h'(\sum_{i=1}^{l-1} \widetilde{e}_{i})$, $b=-h'(-\sum_{i=1}^{l-1}
\widetilde{e}_{i})$, $c=k'(\sum_{i=1}^{l-1} \widetilde{e}_{i})$ and
$d=-k'(-\sum_{i=1}^{l-1} \widetilde{e}_{i})$.  We know the following
inequalities: i) $[x]\leq x\leq [x]+t$, $[y]\leq y\leq [y]+t$ and
$0\leq r\leq t$; ii) $a\leq x\leq b$ and $c\leq y\leq d$; iii)
$a+c\leq [x]+[y]+t=x+y\leq b+d$. Observe that
$\widetilde{m}_{1}(\sum_{i=1}^{l-1} \widetilde{e}_{i})=x$, $m_{1}(
\sum_{i=1}^{l-1} \widetilde{e}_{i})= [x]+r$, $\widetilde{m}_{2}(
\sum_{i=1}^{l-1} \widetilde{e}_{i})= y$ and $m_{2}( \sum_{i=1}^{l-1}
\widetilde{e}_{i})= [y]+t-r$. It is sufficient to show that there is
$r$ such that $0\leq r\leq t$, $a\leq [x]+r\leq b$ and $c\leq
[y]+t-r\leq d$. Observe that  $r$ takes all the value between  0 and
$t$ when  $\bar{r}$ varies between 0 and $l-1$.

1) If $t+[x]\leq b$ we define  $r$ as $\textrm{min}\{[y]+t-c,t\}$. If  $[y]\geq c$ then  $r=t$, so $b\geq [x]+t=
[x]+r\geq x\geq a$ and  $c\leq [y]\leq y\leq d$.  If $c\geq [y]$ then  $b\geq [x]+t\geq [x]+r= [x]+[y]+t-c\geq
a+c-c=a$ and $c=[y]+t-([y]+t-c)=[y]+t-r\leq d$.

2) Suppose now that $[y]+[x]+t\leq b+c$. If $c-[y]$ is positive then we define $r=t+[y]-c$, so $t-r=c-[y]$
($t+[y]\geq y\geq c$ so $r\geq0$). In this case $c=[y]+t-r\leq d$ and $a\leq [x]+[y]+t-c=[x]+r\leq b$. If
$c-[y]$ is negative then we define  $r=t$, so $c\leq [y]=[y]+t-r\leq d$ and $a\leq x\leq [x]+t=[x]+r\leq
c+b-[y]\leq b$.

3) Finally suppose that  $t+[x]> b$ and $[y]+[x]+t>b+c$. We define  $r=b-[x]$, so $a\leq [x]+r=b$ and $d\geq
[y]+[x]+t-b=[y]+t-r\geq c$. $\square$

\

We have proved that the "minimal" weights of $\prod(Z^{n},h+k)$ come
from semi-invariant sections that do not vanish on a suitable
divisor. Moreover we do not need to assume that $L_{h}$ and $L_{k}$
are ample. This fact are no longer true if we consider varieties
whose fans are a little less  symmetric. Notice that the fans of the
varieties considered in the previous proposition are invariant for
any automorphism of $N$ which permutes the vectors of the basis,
fixing $e_{l}$. In the following we define a class of varieties
without such symmetry and obtained by blow-ups from varieties of the
previous family.

\begin{thm}\label{ex 3-2}
Let $\widetilde{Z}^{n}$ be the blow-up of $Z^{n}$  along the stable subvariety associated to
$\sigma(\sum_{j=1}^{l}e_{j},e_{2},...,e_{l})$. Let $L_{h}$ be any ample line bundles on $\widetilde{Z}^{n}$,
then the product of sections $m_{h,h}$ is surjective.
\end{thm}

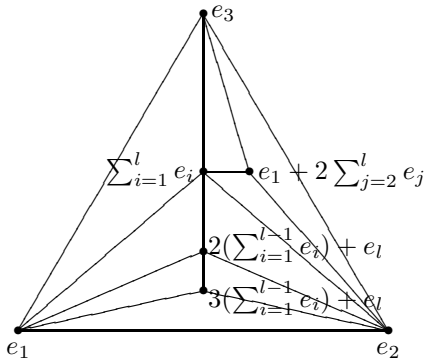
\[\begin{xy} <5pt, 0pt>: (0,0)*{\scriptstyle\bullet}="e1";
(28,0)*{\scriptstyle\bullet}="e2";(14,24)*{\scriptstyle\bullet}="e3"; (14,12)*{\scriptstyle\bullet}="v1";
(14,6)*{\scriptstyle\bullet}="v2"; (14,3)*{\scriptstyle\bullet}="v3";(17.5,12)*{\scriptstyle\bullet}="w";
"e1"+0;"e2"+0**\dir{-};"e1"+0;"e3"+0**\dir{-};"e2"+0;"e3"+0**\dir{-};
"v1"+0;"e1"+0**\dir{-};"v1"+0;"e2"+0**\dir{-};"v1"+0;"e3"+0**\dir{-};
"v2"+0;"e1"+0**\dir{-};"v2"+0;"e2"+0**\dir{-}; "v3"+0;"e1"+0**\dir{-};"v3"+0;"e2"+0**\dir{-};
"v2"+0;"v1"+0**\dir{-};"v2"+0;"v3"+0**\dir{-};
"w"+0;"v1"+0**\dir{-};"w"+0;"e2"+0**\dir{-};"w"+0;"e3"+0**\dir{-};
"e1"+(0,-1.5)*{e_{1}};"e2"+(0,-1.5)*{e_{2}};"e3"+(1.5,0)*{e_{3}}; "v1"+(-4,0.2)*{\sum_{i=1}^{l}e_{i}};
"v2"+(7,0.5)*{2(\sum_{i=1}^{l-1}e_{i})+e_{l}}; "v3"+(7,-0.5)*{3(\sum_{i=1}^{l-1}e_{i})+e_{l}};
"w"+(7,0)*{e_{1}+2\sum_{j=2}^{l}e_{j}};
\end{xy}\]

{\em Proof.} In the figure we have drawn the case in which $l=3$ and $n=3$. We introduce some notation to
simplify the counts: $w:=e_{1}+2\sum_{j=2}^{l}e_{j}$ and $v_{i}:=i(\sum_{i=1}^{l-1}e_{i})+e_{l}$ for each $i$.
Moreover we  suppose that $h(e_{j})=0$ for each $j$. In the proof we allow $L_{h}$ to be the pullback of an
ample linearized line bundle on $Z^{n}$, i.e. $h(w)=h(v_{1})$. We want to prove the proposition by induction on
$h(w)$ and on the dimension of $\widetilde{Z}^{n}$. Observe that if $h(w)= h(v_{1})$ then $m_{h,h}$ is
surjective because of the previous proposition, while if the dimension of $\widetilde{Z}^{n}$ is 2, then
$m_{h,h}$ is surjective because of the theorem~\ref{ex 2}. Suppose now that $h(w)> h(v_{1})$. We want to prove
that $Q_{2h}\cap M=(Q_{h}\cap M)+(Q_{h}\cap M)$ in a similar way to the previous proposition. Let
$a_{i}=h(v_{i})$ and $b=h(w)$. As before, if $p$ belongs to $(Q_{h}\cap M)+(Q_{h}\cap M)$ then
$p+\bigoplus\textbf{Z}^{+}f_{i}$ is contained in $(Q_{h}\cap M)+(Q_{h}\cap M)$. Thus we can suppose that
$m-f_{l}$ does not belong to $(Q_{2h}\cap M)$, so either there is an $i$ such that $m(v_{i})=a_{i}$ or
$m(w)-2h(w)\in \{0,1\}$.

Either we have to study a divisor or $m(w)=2h(w)+1$. One can easily
verify that the only divisor that we do not have already considered
is $\widetilde{Z}^{n}_{\sigma(v_{1})}$. This variety can be studied
in a very similar way to the divisor $Z^{n}_{1}$ of $Z^{n}$, so we
left the details to the reader. Thus we can suppose that
$m(w)=2b+1$. We now want to write some necessary conditions to the
strictly convexity of $h$ on the fan $\Delta$ associated to $Z$. The
conditions $(h|\sigma(v_{1}, w,e_{2},...,e_{l-1})) (v_{i})>
h(v_{i})$, $h|\sigma(w,e_{2},...,e_{l}) (v_{1})>h(v_{1})$,
$h|\sigma(v_{1},e_{1},e_{3},...,e_{l})$ $(e_{2})$ $>h(e_{2})$,
$h|\sigma(v_{1},e_{1},e_{3},...,e_{l})$ $(w)
>h(w)$ and  $h|\sigma(v_{1},e_{1},e_{3},...,e_{l}) (v_{i})> h(v_{i})$ imply:
\[ a_{i}+(i-1)b< (2i-1)a_{1} \ \ \forall i \]
\[ b>a_{1}>0, \ \ \ 2a_{1}>b ,\ \ ia_{1}>a_{i}.\]

Let $\Delta$ be the fan of $\widetilde{Z}^{n}$, let $\Delta'$ be the fan of $Z^{n}$ and let $h'$ be the $\Delta$
linear function  such that $h'(e_{i})=0$, $h'(v_{i})= h(v_{i})$ and $h'(w)= h(w)-1$. We need the following easy
lemma on $h'$ which we will not prove.

\begin{lem}
$h'$ is convex on $\Delta$ and is strictly convex either on  $\Delta$ or on $\Delta'$.
\end{lem}

By induction we can suppose that $m_{h',h'}$ is surjective, so we can suppose that there are two points
$m_{1}\in Q_{h}\cap M$ and $m_{2}\in Q_{h'}\cap M$ such that $m_{1}+m_{2}=m$. If $m_{2}$ belong to $Q_{h}$ then
$m\in Q_{h}\cap M+Q_{h}\cap M$. Otherwise we have $m_{2}(w)=b-1$ and $m_{1}(w)=b+2$. Write
$m_{1}=(x_{1},...,x_{l})$ and $m_{2}=(y_{1},...,y_{l})$.

We can suppose that $m_{1}-f_{l}\ \nin\ Q_{h}$ because  $m_{2}+f_{l}\in Q_{h}$. Thus there is  $i$ such that
$m_{1}(v_{i})=a_{i}$. Moreover we can suppose that $(m_{1}+f_{1}-f_{j}, m_{2}-f_{1}+f_{j})$ does not belong to
$Q_{h}\times Q_{h'}$ for any $j=2,...,l-1$, so  $x_{j}=0$ or $y_{1}=0$. If $y_{1}=0$ then $2a_{1}-1\geq
b-1=m_{2}(w)=2\sum y_{j}=2m_{2}(v_{1})\geq 2a_{1}$, so we have obtained a contradiction. Hence $y_{1}\neq0$ and
$x_{j}=0$ for each $j=2,...,l-1$. Suppose that there is $i>1$ such that $m_{1}(v_{i})=ix_{1}+x_{l}=a_{i}$, then
we have

\[ (2i-1)a_{1}\leq
 (2i-1)(x_{1}+x_{l})=m_{1}(v_{i}) +(i-1)m_{1}(w)=\]
\[=a_{i}+(i-1)(b+2)\leq
 (2i-1)a_{1}+2i-2,\]
so $0\leq (2i-1)(x_{1}+x_{l}-a_{1})\leq 2i-2 $ (remember that $a_{i}+(i-1)b<(2i-1)a_{1}$). We have
$x_{1}+x_{l}=a_{1}$ because $x_{1}+x_{l}-a_{1}$ is an integer. Observe that we have showed that
$m_{1}(v_{1})=x_{1}+x_{l}=a_{1}$ or $m_{1}(e_{l})=x_{l}=0$. In the last case we have $x_{2}=...=x_{l}=0$ and
$x_{1}=b+2$. We can suppose that $(m_{1}-f_{1},m_{2}+f_{1})$ does not belong to $Q_{h}\times Q_{h'}$, so there
is $s>0$ such that $m_{1}(v_{s})-a_{s}<s$. Observe that $m_{1}(v_{s})=sx_{1}=sb+2s$, so $a_{s}\leq
sb=m_{1}(v_{s})-2s<a_{s}-s<a_{s}$, so we have obtained a contradiction.

Finally we can suppose that $x_{j}=0$ for each  $j=2,...,l-1$, $x_{1}+x_{l}=a_{1}$ and $x_{1}+2x_{l}=b+2$, so
$x_{1}=2a_{1}-b-2$ and $x_{l}=b+2-a_{1}$. Moreover we can suppose that $(m_{1}+f_{1}-f_{l}, m_{2}-f_{1}+f_{l})$
does not belong to $Q_{h}\times Q_{h'}$, so $ x_{l}=0$, $y_{1}=0$ or there is  $i>1$ such that
$\varepsilon:=m_{2}(v_{i})- a_{i}<i$. Observe that we have already considered the first two cases.

We have $a_{i}\leq m_{1}(v_{i})= ix_{1}+x_{l}= (2i-1)a_{1}-(i-1)b-2(i-1)$, so $(2i-1)a_{1}\geq
a_{i}+(i-1)b+2(i-1)$. Finally

\[ (2i-1)a_{1}\leq (2i-1)(\sum y_{j})\leq
 (2i-1)y_{1}+ (3i-2)\sum_{j\neq 1,l} y_{j}+(2i-1)y_{l}=\]
\[= m_{2}(v_{i})+(i-1)m_{2}(w)
=(i-1)(b-1)+a_{i}+\varepsilon=\]
\[=(i-1)b+a_{i}+2(i-1)-3(i-1)+\varepsilon\leq (2i-1)a_{1}-3(i-1)+\varepsilon,\]
so $3(i-1)\leq\varepsilon\leq i-1$, a contradiction. $\square$


\end{document}